\documentclass[12pt,twoside]{article}
\usepackage{amsfonts,a4,epsf}

\newcommand{\frakg}{\ensuremath\mathfrak g}

\newcommand{\frakf}{\ensuremath\mathfrak f}
\newcommand{\frake}{\ensuremath\mathfrak e}
\newcommand{\fraksl}{\ensuremath{\mathfrak {sl}}}
\newcommand{\frakso}{\ensuremath{\mathfrak {so}}}
\newcommand{\fraksp}{\ensuremath{\mathfrak {sp}}}
\newtheorem{prop}{Proposition}[section]

\newtheorem{ex}[prop]{Example}
\newtheorem{lem}[prop]{Lemma}
\newtheorem{rmk}[prop]{Remark}
\newtheorem{df}[prop]{Definition}
\newtheorem{dt}[prop]{Theorem and Definition}

\newcommand{\gf}{{\mathfrak g}}

\newcommand{\C}{{\bf C}}
\newcommand{\R}{{\bf{R}}}
\newcommand{\topf}{\ensuremath {{\cal T}}_{tpf} }
\newcommand{\topfc}{\ensuremath{{\cal T}}_{topfc} }

\newcommand{\I}{I_{(\mathfrak g,V)}}
\newcommand{\Ig}{I_{({\mathfrak g}_2,V)}}
\newcommand{\If}{I_{({\mathfrak f}_4,V)}}
\newcommand{\Ie}{I_{({\mathfrak e}_6,V)}}
\newcommand{\Iee}{I_{({\mathfrak e}_7,V)}}
\newcommand{\Ieee}{I_{({\mathfrak e}_8,L)}}

\newcommand{\openbox}{\leavevmode
  \hbox to.77778em{%
  \hfil\vrule
  \vbox to.675em{\hrule width.6em\vfil\hrule}%
  \vrule\hfil}}
\newcommand{\schluss}{\vspace{-4ex}  \hfill \ensuremath{\openbox} \vspace{4ex}}
\renewcommand{\Box}{\schluss}

\newcommand{\Col}{{\cal C}ol}
\newcommand{\col}{\ensuremath col}
\newcommand{\Dc}{{\cal D}_{3c}}
\newcommand{\ra}{\rightarrow}
\newcommand{\la}{\leftarrow}

\newcommand{\pg}[3]{p^\frakg_{#1\otimes #2,#3}}
\newcommand{\ig}[3]{i^\frakg_{#1,#2\otimes #3}}
\def\objez{$\begin{array}{l}
c_1^{x_1}\\[1ex]
c_2^{x_2}
\end{array}$}
\def\objzd{$\begin{array}{l}
c_2^{x_2}\\[1ex]
c_3^{x_3}
\end{array}$}
\def\objze{$\begin{array}{l}
c_2^{x_2}\\[1ex]
c_1^{x_1}
\end{array}$}
\def\objx{$\begin{array}{l}
c^{\overline{x}}\\[1ex]
c^{x}
\end{array}$}

\newcommand{\pstext}[1]{\hspace{1mm}\raisebox{-0.5ex}{\epsfysize=2.5ex
\epsffile{#1.eps}}\hspace{1mm}}
\newcommand{\psdiag}[3]{\hspace{1mm}\raisebox{-#1mm}{\epsfysize#2mm
\epsffile{#3.eps}}\hspace{1mm}}

\title{Skein relations for the link invariants coming from exceptional Lie algebras}
\author{Anna-Barbara Berger\footnote{The authors are partially supported by the
Schweizerische Nationalfonds.}\hspace{0.2em} and Ines Stassen$^*$}
\date{September 1998}

\begin{document}
\maketitle

\begin{abstract}
Pulling back the weight systems associated with the exceptional
Lie algebras and their standard representations by a modification
of the universal Vassiliev-Kontsevich invariant yields link
invariants; extending them to coloured 3-nets, we derive for each
of them a skein relation.
\end{abstract}

\setcounter{section}{-1}
\section{Introduction}

There is a well-known technique for the construction of Vassiliev link
invariants: define a weight system (i.e.\ a linear form on the space of chord
diagrams respecting certain relations)
on the basis of some Lie algebraic data and pull it back by
the universal Vassiliev-Kontsevich invariant. But unfortunately, the latter
is not known explicitly enough to allow direct evaluation of these link
invariants.

Efforts have been made to handle the universal Vassiliev-Kontsevich
invariant by considering only
``elementary'' parts of links into which any link may be cut.
This approach has been successful in so far as one
may hope to find skein relations for the link invariants coming from Lie
algebras---a skein relation is an equation implying a recursive algorithm
for the computation of a link invariant, for example the following
equation, which determines the famous Jones polynomial up to normalization:
$$t^2 P(\psdiag{2}{6}{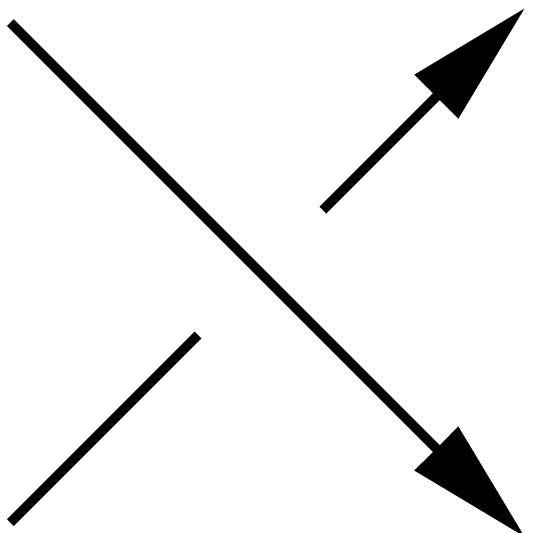})-t^{-2}P(\psdiag{2}{6}{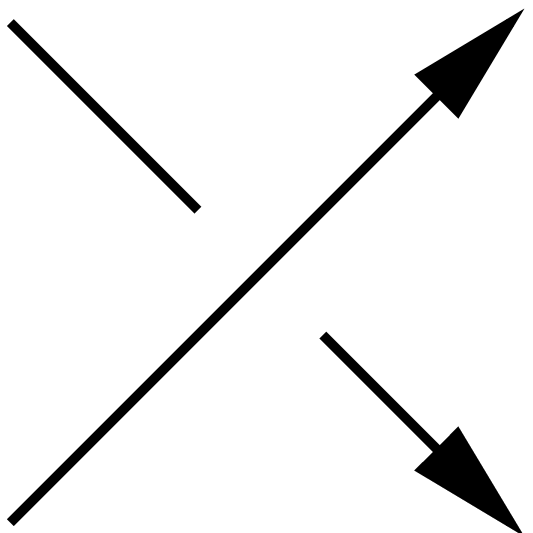})
= (t^{-1}-t)P(\psdiag{2}{6}{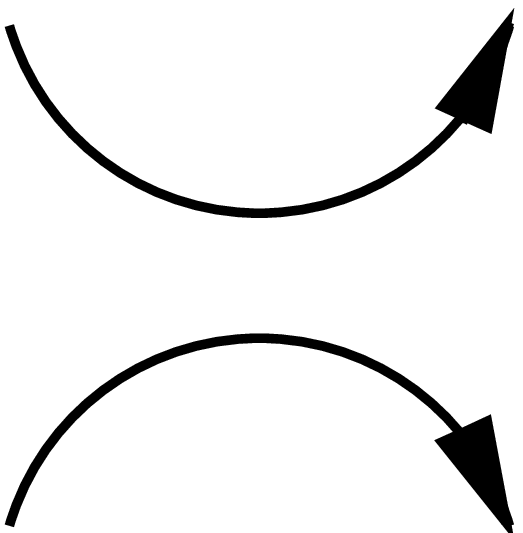}).$$
\medskip

It has been shown that the link invariants obtained from
the classical simple Lie algebras ${\fraksl}_n$, ${\frakso}_n$, and ${\fraksp}_n$
satisfy certain
versions of the skein relation of the HOMFLY polynomial (${\fraksl}_n$; see [LM 1])
resp.\ the Kauffman polynomial (${\frakso}_n$, ${\fraksp}_n$; see [LM 2]).
But what about the exceptional simple Lie algebras?
\medskip

In [BS], we have outlined a strategy for establishing skein
relations and dealt with the case of the exceptional Lie algebra
$\mathfrak g _2$. At the price of further generalizations of the
notion of links---besides branchings, we have to introduce a
colouring---we can now present skein relations for all the
invariants coming from exceptional simple Lie algebras and their
standard representations.
Unfortunately, we have not been able to determine all the initial
values of the recursive algorithms implied by the skein relations.
\medskip

To the reader who is not familiar with Lie theory, we recommend [H]
and [FH]. For an introduction to Vassiliev invariants and weight systems,
see [BN 1]; a more general definition of weight systems is given in [V],
section 6.
\medskip

{\bf Overview} over the categories and functors appearing in this paper:
\vspace{-2mm}\\
\begin{center}
$\I$ \\
\psdiag{0}{14}{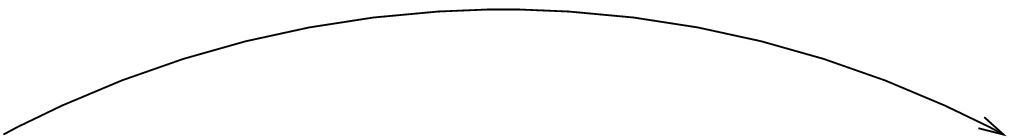} \vspace{-1.2cm}\\
\footnotesize the $(\frakg,V)$-invariant
\end{center}
\vspace{0mm}
\begin{tabular}{@{} c c c c c}
\hspace{2.85cm} & \hspace{2cm} & \hspace{2.9cm} & \hspace{2.2cm}
& \hspace{2.8cm} \\
${\cal T}_{topfc}$ & $\stackrel{\hat Z_f}{\longrightarrow}$ & $\hat{\cal D}_{3c}$ &
$\stackrel{\Psi_\frakg}{\longrightarrow}$ & ${\cal C}(\frakg)$
\end{tabular}
\vspace{0mm}\\
\begin{tabular}{@{} p{2.85cm} p{2cm} p{2.9cm} p{2.2cm} p{2.8cm}}
\footnotesize the ca\-te\-go\-ry of {\bf t}ri\-va\-lent, {\bf o}rien\-ted,
{\bf p}a\-ran\-the\-sized,
{\bf f}ramed, {\bf c}o\-loured tangles  &
\footnotesize the universal Vassi\-liev-Kontsevich in\-va\-riant &
\footnotesize the ca\-te\-go\-ry of co\-loured
3-dia\-grams (a ge\-ne\-ra\-li\-za\-tion of chord dia\-grams) &
\footnotesize the $(\frakg,V)$-weight system &
\footnotesize a sub\-ca\-te\-go\-ry of the ca\-te\-go\-ry of
re\-pre\-sen\-ta\-tions of $\frakg$
\end{tabular}
\vspace{6mm}

{\bf Acknowledgements} \hspace{0.3em} We would
like to express our thanks to Ch.\ Riedtmann
for supervision and encouragement and to P. Vogel
for suggesting to work on this subject.

\section{Invariants for coloured oriented 3-tangles}
In this section, we generalize the terms 3-nets, 3-tangles, and
3-diagrams introduced in [BS] by adding orientations and colourings, and we adapt the universal
Vassiliev-Kontsevich and the $(\frakg,V)$-weight system to these generalizations.
This will allow us to obtain for each exceptional Lie algebra a skein relation
for a corresponding link invariant.\\

A 3-net is something like a ``link with branchings''. To describe
the situation near a trivalent vertex (i.e.\ near a branching
point), we will need the following notion:\\ Let $B$ be the open unit
ball in $\R^3$, i.e.\ $\{x\in
\R^3; |x|< 1\}$, together with the distinguished subset $T:= \{(t,0,0)|0\leq t\leq1\}\cup
\{(-\frac{1}{2},\frac{\sqrt{3}}{2},0)t|0\leq t\leq 1\}\cup
\{(-\frac{1}{2},-\frac{\sqrt{3}}{2},0)t|0\leq t\leq 1\}$.\\[1em]
\centerline{$B\quad =\quad$\psdiag{7}{21}{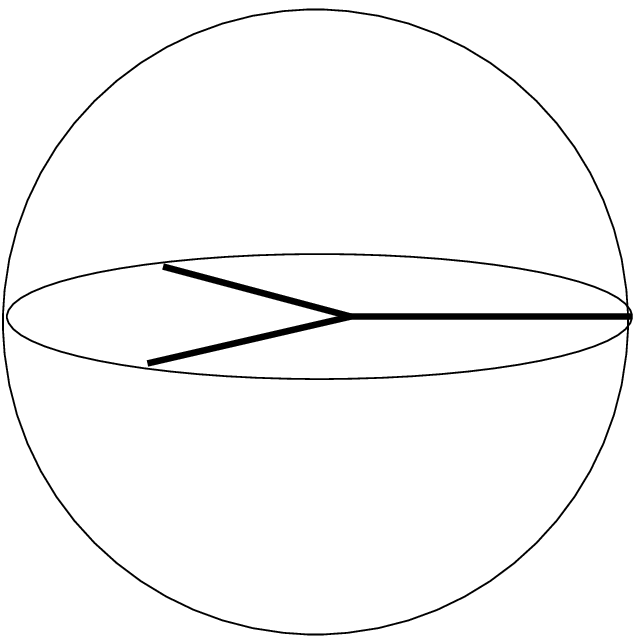}}

\begin{df}
An {\em oriented framed 3-net} is a subset $N$ of \/ $\R^3$ with a
finite subset $\{t_1,\ldots ,t_n\}\subset N$ such that:
\begin{itemize}
\item[(i)] there exist disjoint open subsets $U_1,\ldots,U_n$ and
diffeomorphisms $f_i:U_i\rightarrow B\,(i=1,\ldots,n)$ such that
$U_i$ is a neighbourhood of $t_i$, $f_i(t_i)=(0,0,0)$, and $f_i(N\cap U_i)= T$,
\item[(ii)] $\tilde{N}:= N \backslash \left(\bigcup_{i=1}^n f_i^{-1}(\{x\in B;
|x|<\frac{1}{2}\})\right)$ is an embedded  smooth  closed compact 1-dimensional
manifold,
\end{itemize}
together with:
\begin{itemize}
\item[(iii)] an orientation on $N \backslash \{t_1,\ldots,t_n\}$,
\item[(iv)] a  smooth vector field on $N$ that is nowhere tangent to $N$
(and in particular nowhere zero)
.\\
\end{itemize}

  The points
$t_1,\ldots, t_n$ are called {\em trivalent vertices} of $N$;
boundary points $x$ of $\tilde{N}$ with $x\not\in U_i
\:(\forall i)$ are called {\em univalent vertices} of $N$. An {\em
edge} of $N$ is a connected component of $N \backslash
\{t_1,\ldots,t_n\}$.
\end{df}
\begin{df}
Let $\Col :=\{so(lid), da(shed),do(tted)\}$ be the set of colours.
\end{df}
\begin{df}
An {\em oriented framed coloured 3-net} is an oriented framed 3-net
together with a colouring of the edges, i.e. a map $\col:\{
\mbox{edges of } N \} \rightarrow \Col.$
\end{df}

The definitions for {\it equivalent, closed}, and {\it planar}
oriented framed coloured 3-nets and 3-tangles are the adapted
versions of the corresponding definitions for framed 3-nets and
3-tangles in [BS].

The category $\topfc$ of trivalent oriented parenthesized framed
coloured tangles is the coloured and oriented analogon of the
category $\topf$ in [BS]. Its objects are non-associative coloured
words:
\begin{df}
A {\em non-associative coloured word} is a word $w$ in the alphabet
$\{),(\}\cup\{c^\rightarrow,c^\leftarrow| c \in
\Col\}$, such that $w$ is equal to the empty word or to $(c^\ra),( c^\la)$,
 or $(w_1w_2)$
where $c$ is any colour and $w_1, w_2$ are non-associative coloured words.
For every non-associative coloured word $w$, we identify $(w)$ with $w$.\newline For $x
\in \{\rightarrow,\leftarrow\}$ let $\overline{x}$ be $\leftarrow$
if $x=\,\rightarrow$ and $\rightarrow$ if $x=\,\leftarrow$.\newline The
{\em underlying string $u(w)$ of a non-associative coloured word}
$w$ is the sequence in $c^\rightarrow, c^\leftarrow$ with $c\in
\Col$ one obtains from $w$ by omitting all parentheses.\newline The
{\em length} $l(w)$ of a non-associative coloured word is the
number of symbols in $u(w)$.\end{df}
\begin{ex}{\rm
$\hspace*{1.7cm}u\left((((c_1^\ra)( c_2^\la))(((c_3^\la)( c_2^\ra))(c_1^\la)))\right)=c_1^\ra
c_2^\la c_3^\la c_2^\ra c_1^\la$
$$l\left((((c_1^\ra)(c_2^\la))(((c_3^\la)(
c_2^\ra))(c_1^\la)))\right)=5
$$
}\end{ex}

\begin{df}
Let $\topfc$ \/ be the monoidal $\C$-category that is given by the
following data:\newline {\bf objects:} non-associative coloured
words. The unit object is the empty word, and the tensor product on
the objects is defined by $w_1\otimes w_2:=(w_1w_2)$.\newline {\bf
morphisms:}\\ {\bf generators:} The morphism spaces are generated
by:
\begin{itemize}
\item[(G1)] A morphism ${\mbox{\pstext{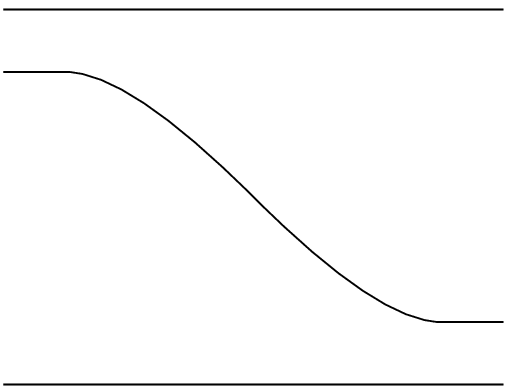}}}_{v,w,x}$ and a morphism
 ${\mbox{\pstext{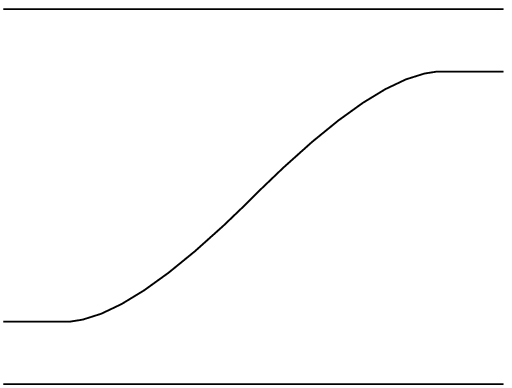}}}_{v,w,x}$ for each triple $(v,w,x)$ of
non-empty non-associative coloured words. The sources of these
morphisms are $((vw)x)$ and $(v(wx))$ and their targets are $(v(wx))$
and $((vw)x)$ respectively.
\item[(G2)]
A morphism ${\mbox{\pstext{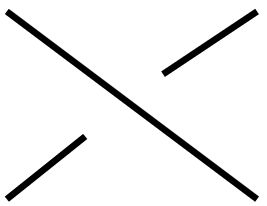}}}_{v,w}$ and a morphism
 ${\mbox{\pstext{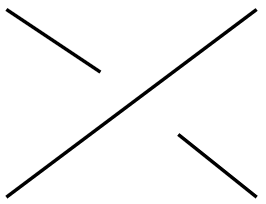}}_{v,w}}$ for each pair $(v,w)$ of
non-empty non-associative coloured words. The source of these morphisms is
$(vw)$ and their target is $(wv)$.
\item[(G3)] A morphism ${\mbox{\pstext{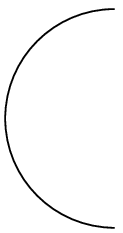}}}_c$ and a morphism $\pstext{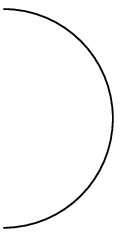}_c$
with source {\it empty word} and $((c^{\overline{x}})(c^x))$ and target $((c^{\overline{x}})(c^x)$
and {\it empty word} respectively for each $c\in \Col$ and $x\in
\{\ra,\la\}$.
\item[(G4)] A morphism ${\mbox{\pstext{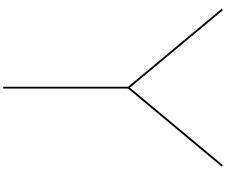}}}_{c_1^{x_1},c_2^{x_2},c_3^{x_3}}$
 and a morphism ${\mbox{\pstext{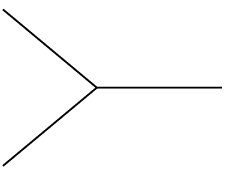}}}_{c_1^{x_1},c_2^{x_2},c_3^{x_3}}$ with source
$(c_1^{x_1})$ and $((c_1^{x_1})( c_2^{x_2}))$ and target $((c_2^{x_2})(c_3^{x_3}))$ and $(c_3^{x_3})$
respectively for each triple $(c_1^{x_1},c_2^{x_2},c_3^{x_3})$
with $c_i \in \Col$ and $x_i \in \{\rightarrow,\leftarrow\}.$\\
\end{itemize}
{\bf relations:}\newline To every morphism of $\topfc$ we can
assign an oriented framed coloured 3-tangle: Proceed like in the
case of ${{\cal T}}_{tpf}$ (see [BS]) to obtain an uncoloured, unoriented
3-tangle and then colour each edge with the colour indicated by
source and/or target and orient them as indicated in the exponent
of the components of the objects.\\ We impose the following
relation on the morphism spaces: two morphisms from $u$ to $w$ are
equivalent if they get assigned equivalent oriented framed coloured
3-tangles.
\end{df}

\begin{ex}{\rm
$(id_{(((so^\ra)(da^\la))(so^\la))}\otimes
\pstext{trir}_{so^\ra,so^\la,do^\ra})\pstext{Zarg5}_{(so^\ra),(((so^\ra)
(da^\la))(so^\la))}$\\ gets assigned:\hspace{1cm}
\psdiag{6}{18}{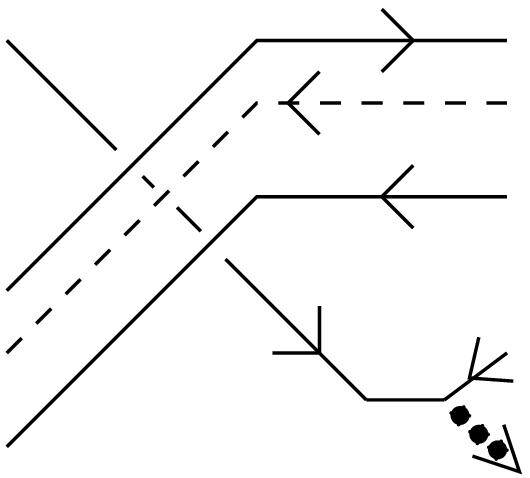}
}\end{ex}

{\bf Convention:} In a graphical representation of a morphism, the
colouring of the edges can either be read from the adjacent
component(s) of their source and/or target or, if these are
omitted, is indicated by the style of the line.\\
 Example: $da^\ra
\psdiag{2}{6}{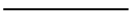}\, da^\ra = \psdiag{0.8}{3}{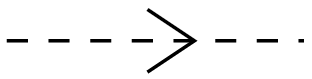}$.

\begin{df}
A {\em coloured oriented 3-diagram} is a finite trivalent graph K (by which we understand
a graph with every
vertex being either univalent or trivalent or else bivalent and adjacent
to a loop) equipped with the following data:
\begin{itemize}
\item for every trivalent vertex x of K, a cyclic order of the edges
arriving at x.
\item a colouring of the edges of $K$,
i.e. a map $\col :\{\mbox{edges of }K\} \rightarrow \Col$.
\end{itemize}

The {\em degree} of a 3-diagram is the number of trivalent vertices
adjacent to at least one dashed edge\footnote{Note that
for a 3-diagram without univalent vertices adjacent to a chord,
this is twice the classical degree.}.
\end{df}

Usually, we describe the 3-diagrams by graphical representations
in the plane encoding
the information about the cyclic order near the trivalent vertices by
arranging the adjacent edges counterclockwise.

\begin{df}
The category $\Dc$ is a monoidal\/ {\bf C}-category whose morphisms are
linear combinations of
certain graphical representations of coloured oriented 3-diagrams. It is given by the
following data:\newline
{\bf objects:} {\em Obj}$(\Dc) := \,\,^{\cdot}\hspace{-2.8mm}\bigcup_{n=0}^\infty\{c^\ra,c^\la | c\in \Col\}^n$.
The tensor product on {\em Obj}$(\Dc)$ is the
juxtaposition.\newline
{\bf morphisms:}\newline
{\bf generators:} The morphism spaces are generated by:
\begin{tabbing}
\objez \hspace{-1mm}\psdiag{3.5}{9}{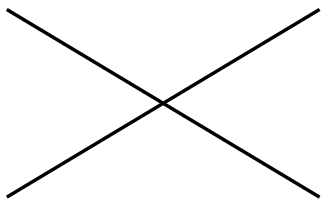}\objze \hspace{2cm}\=\objx\hspace{-1mm}\psdiag{3}{8}{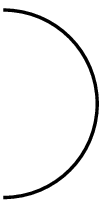}\hspace{2cm}
\=\objez\hspace{-1mm}\psdiag{3}{7}{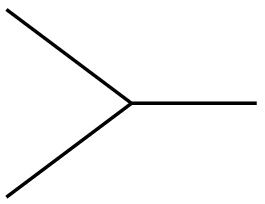}$c_3^{x_3}$\\
\>\psdiag{3}{8}{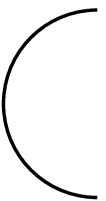}\hspace{-1mm}\objx\>$\,\,\,c_1^{x_1}$\psdiag{3}{7}{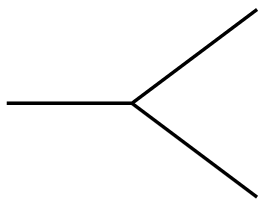}\objzd
\end{tabbing}
the source (resp.\ the target) being denoted on the left-(resp.\
right-)hand side from top to bottom and the edges oriented and
coloured as indicated by the objects.
\footnote{The apparent 4-valent vertices are no vertices at all - they are just
crossings of two edges (there is no need to say that one of them passes
over the other).}
\\ The tensor product of two morphisms is obtained by putting the first above the
second, the composition by glueing together the corresponding entries of
the target of the first and the source of the second.
\\
{\bf relations:} Of course, different graphical representations of
isomorphic co\-lou\-red oriented 3-dia\-grams are to represent the same morphism;
in addition, we impose the following re\-la\-tions:
\begin{tabbing}
(AS1) xxxxxxx \= \kill
(AS) \> $\begin{array}{r}c^x\\ \\  \\ \end{array}\!\psdiag{3}{9}{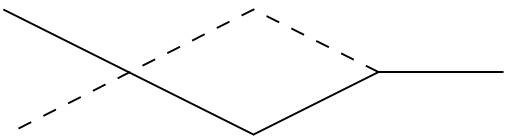}c^x \qquad = \qquad -
\begin{array}{r}c^x\\ \\  \\ \end{array}\!\psdiag{3}{9}{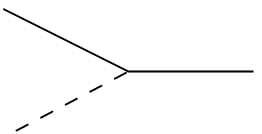}c^x $\\[2ex]
(IHX) \> $\psdiag{4.5}{13.5}{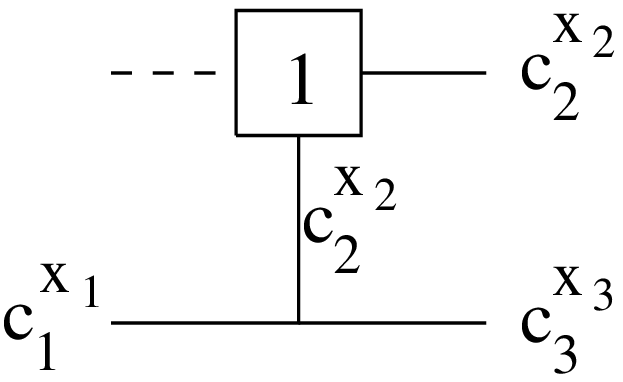}\qquad = \qquad \psdiag{5.5}{16.5}{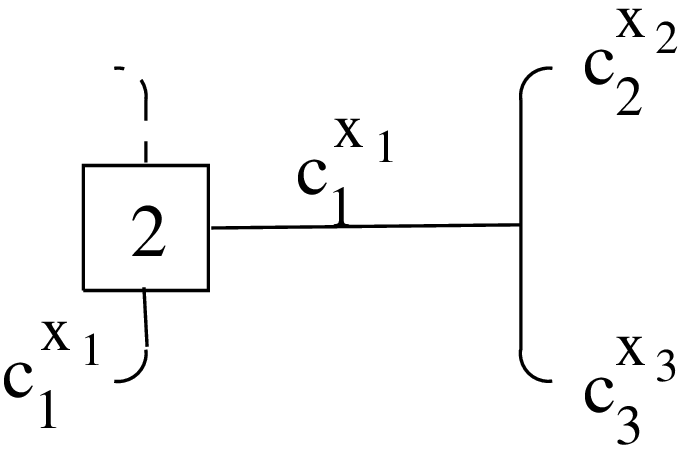}  -  \psdiag{6}{18}{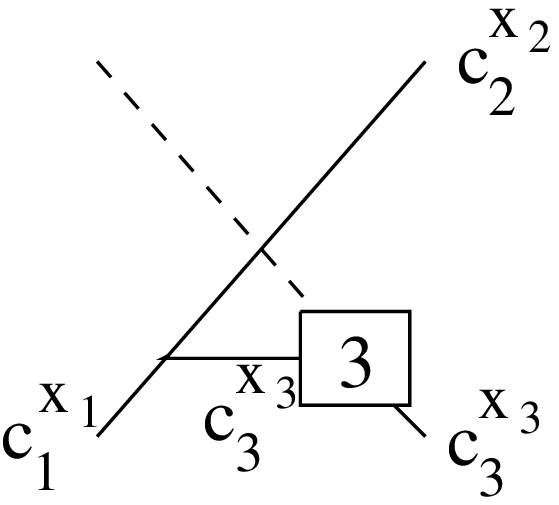}$
\end{tabbing}
\begin {tabbing}
with\hspace{3cm}\= \psdiag{2.5}{7.5}{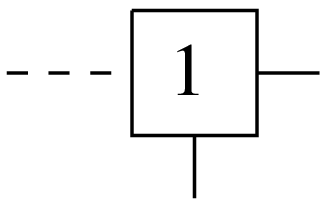}\quad\= = \quad
\=\psdiag{1.8}{5.4}{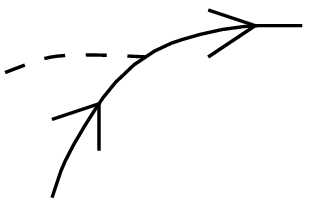}\quad\=rsp.\quad\=\psdiag{2}{6}{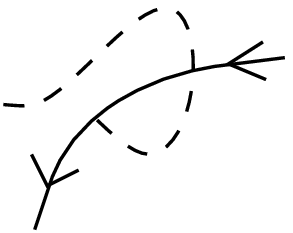}\\[1ex]
\> \psdiag{3.3}{9.9}{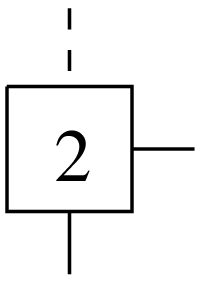}\quad\> = \quad
\>\psdiag{3.4}{10.2}{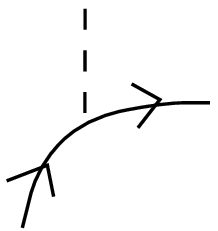}\quad\>rsp.\quad\>\psdiag{3}{9}{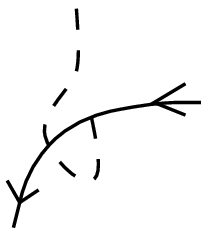}\\[1ex]
\> \psdiag{3}{9}{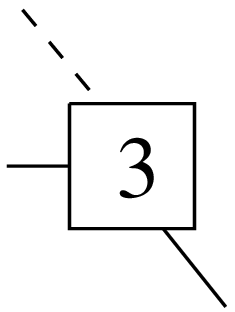}\quad\> = \quad \>\psdiag{3}{9}{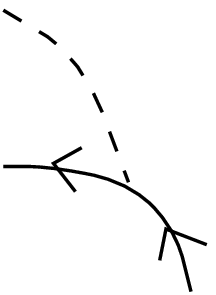}\quad\>rsp.\quad\>\psdiag{3}{9}{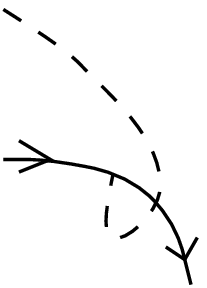}\\
\end{tabbing}
with all possible orientations of the edges and all choices of
colours $c,c_1,c_2,c_3 \in \Col$.
\end{df}

\begin{df}
Let $\hat{\cal D}_{3c}$ be the completion of the (graded) category $\Dc$.
\end{df}

For convenience of notation, we define a family of maps $\Delta_w$.
To do this, we need the following notions:\\

Let a {\em connection} in a coloured oriented 3-diagram be a
sequence of solid edges such that the first and the last edge are
adjacent to a univalent vertex and the head of the first and the
tail of the second of any two subsequent edges is identical and the
third edge arriving at this trivalent vertex is dashed.\newline
A morphism of $\hat \mathcal{D}_{3c}$ may be equipped with some
extra information: in each of its terms, some connections may be
labelled. By this, we obtain what we will call a {\em labelled}
morphism of $\hat \mathcal{D}_{3c}$.\newline
If $D$ is a morphism of $\hat \mathcal{D}_{3c}$ whose source has as
$k$-th component a ``solid'', we denote by $D_k$ the labelled morphism
that is given by $D$ and the labelling of the connection parting from the
$k$-th component of the source in each summand.\newline
For any non-associative coloured word $w$, call $\Delta_w$ the
map which is defined on the space of all labelled morphisms of
$\hat \mathcal{D}_{3c}$ as follows:\\ In every summand, the
labelled connection is replaced by a bunch of parallel strands the
uppermost of which is coloured and oriented according to the first
component of $w$, the second according to the second component of
$w$, and so on; then the sum over all possible ways of lifting the
trivalent vertices of the original connection to the new ones is
taken.

\begin{ex}
\begin{itemize}
\item[]
\item[(a)] $w:=(so^{\to}, da^\leftarrow)$,
$D:=\psdiag{2}{6}{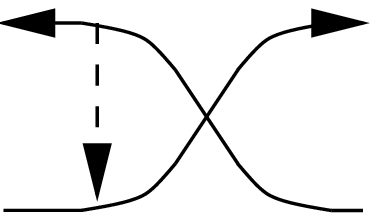}$,
$\Delta_w(D_1) = \psdiag{3}{9}{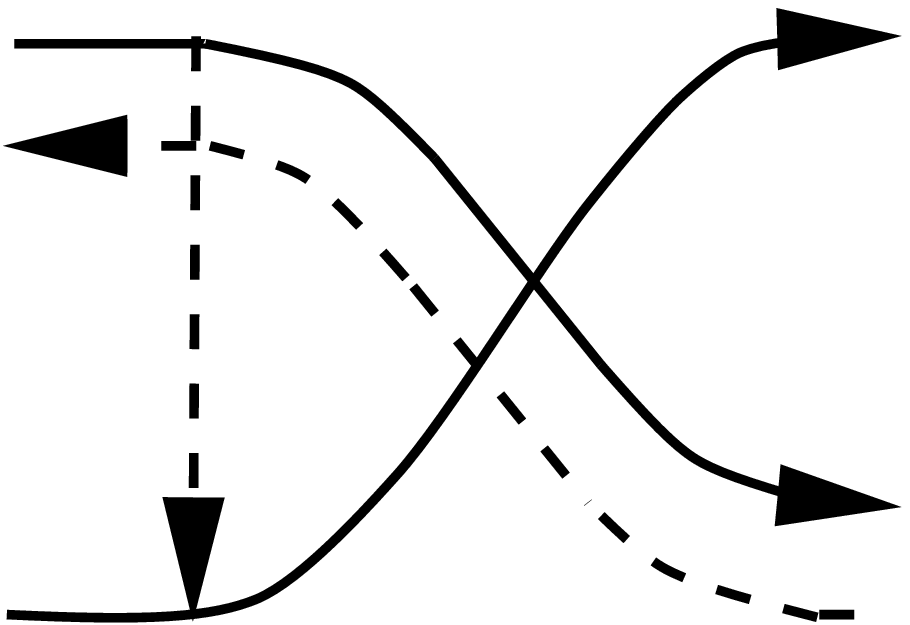} + \psdiag{3}{9}{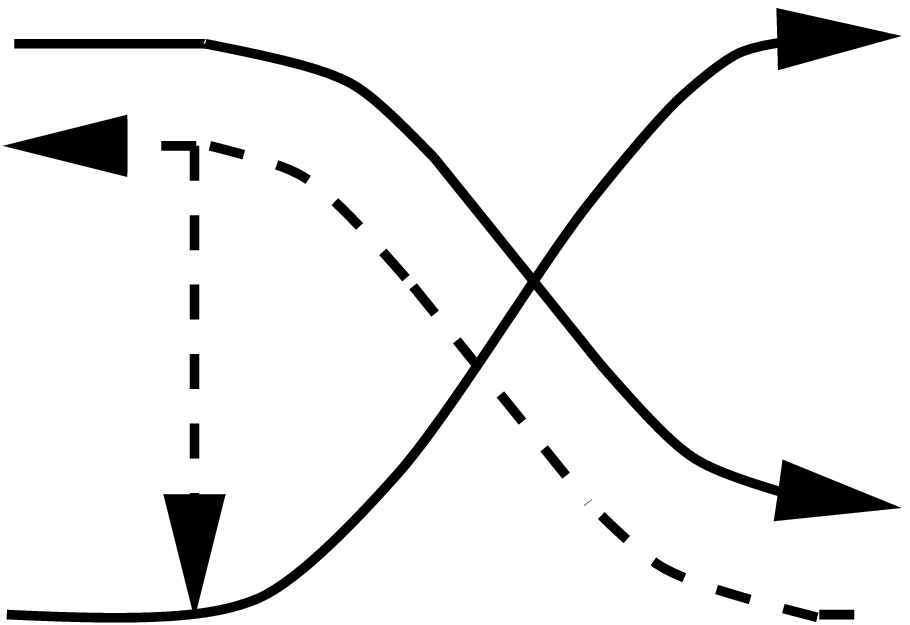}$
\item[(b)] $w:=(so^\to, do^\leftarrow$),
$D:=\psdiag{2}{6}{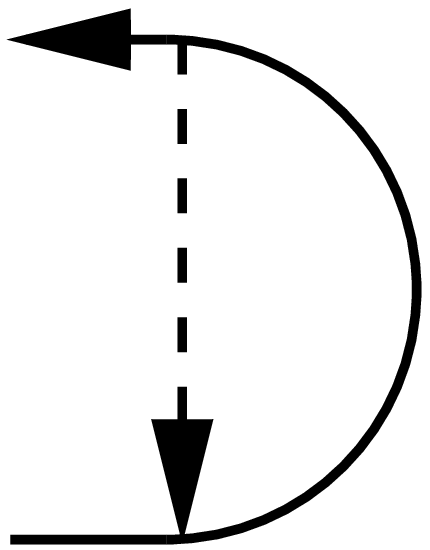}$,
$\Delta_w(D_2) =
\psdiag{3}{9}{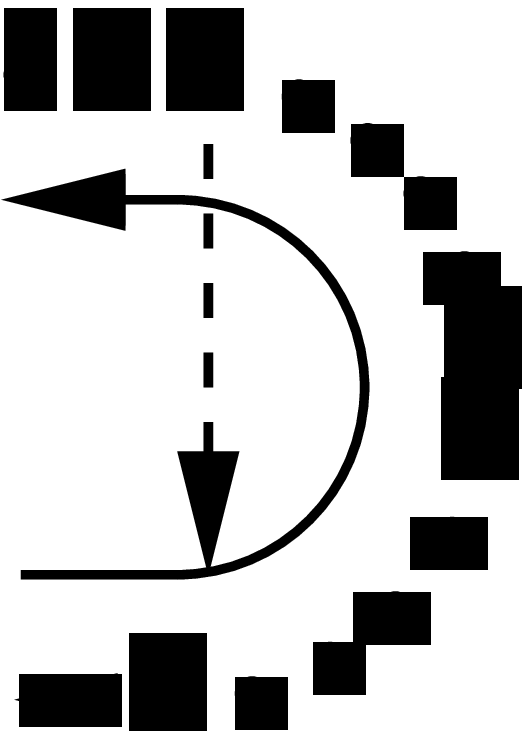} + \psdiag{3}{9}{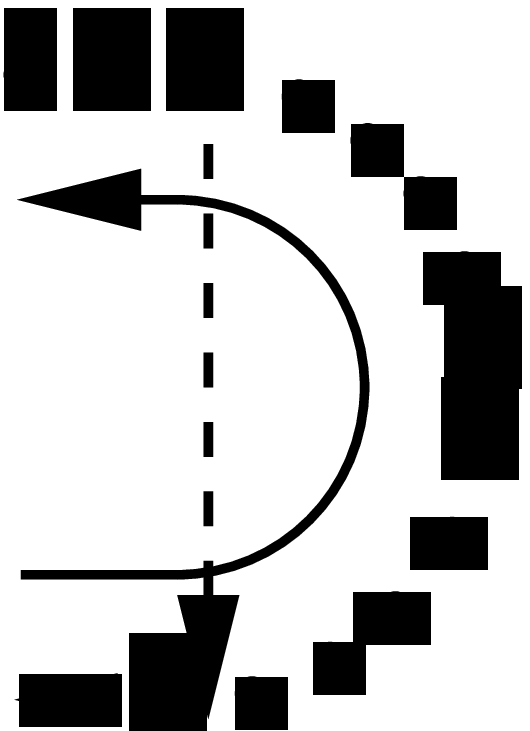}
+ \psdiag{3}{9}{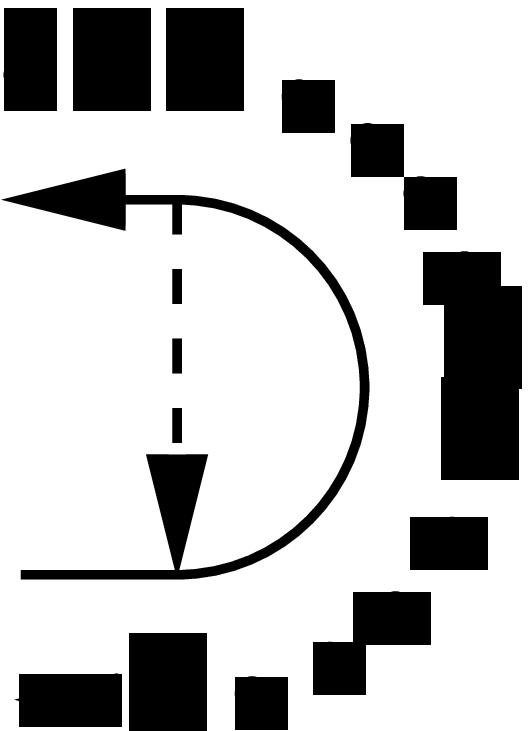}$ + \psdiag{3}{9}{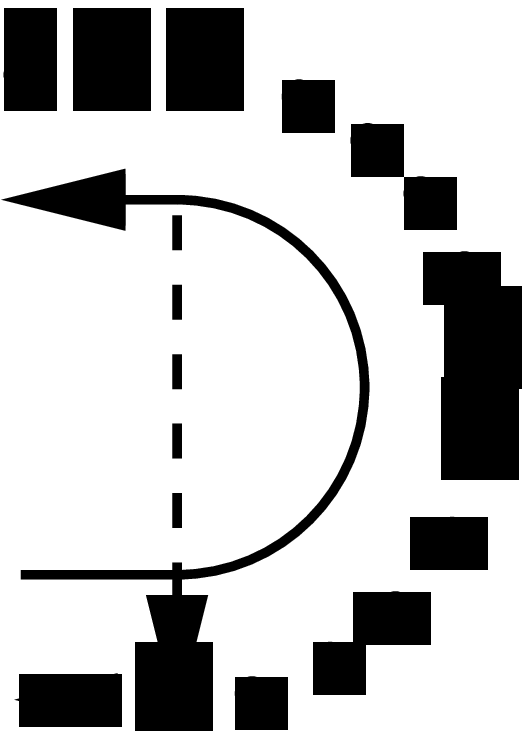}
\end{itemize}
\end{ex}

As an immediate consequence of the relation {\em (IHX)}, we obtain:

\begin{lem}\label{lemma2}
$$\mbox{\objez}\psdiag{9}{15}{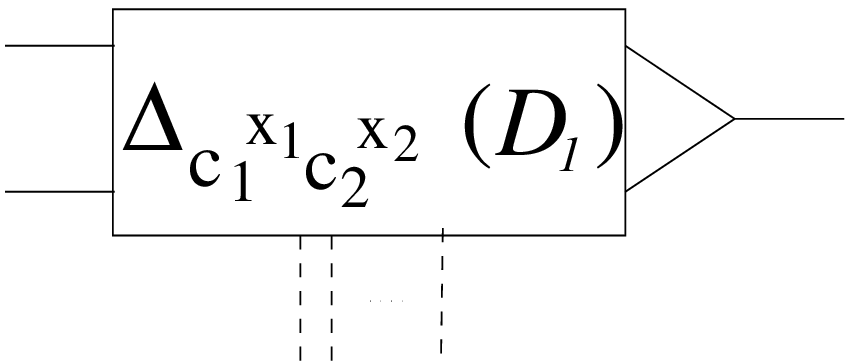}c_3^{x_3} =
\mbox{\objez}\psdiag{10}{15}{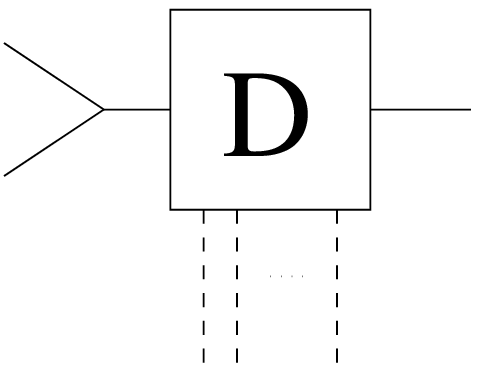}c_3^{x_3}$$
$$c_1^{x_1}\psdiag{9}{15}{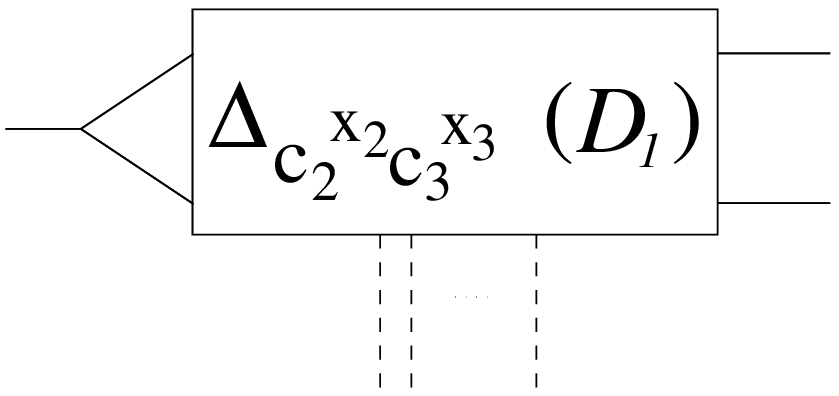}\mbox{\objzd} = \, c_1^{x_1}\psdiag{10}{15}{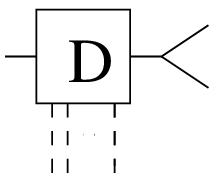}\mbox{\objzd}$$
 for any fitting 3-diagram D. \hfill $\Box$
\end{lem}

Finally, we get to the definition of the universal Vassiliev-Kontsevich
invariant:

\begin{dt}\label{VK}
The following assignments define a monoidal
functor $\hat Z_f : {\cal T}_{topfc} \to \hat{\cal D}_{3c}$, the
{\em extension of the  universal Vassiliev-Kontsevich in\-va\-riant} to
oriented framed coloured 3-tangles:
\vspace{2mm}\\
$\hat Z_f(w):= u(w)$, where $w$ is a non-associative coloured word.
\begin{tabbing}
$\hat Z_f$(\psdiag{3}{9}{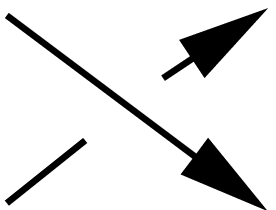})xx \= := xx\=\psdiag{3}{9}{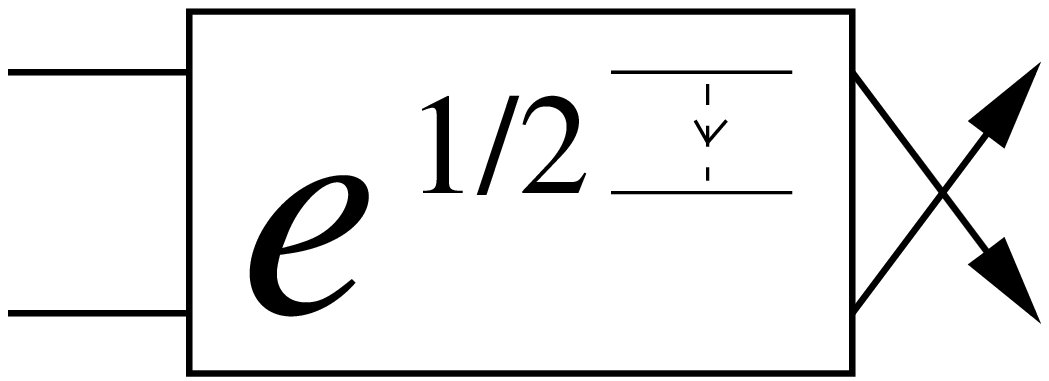} xxxxxx
\= $\hat Z_f$(\psdiag{3}{9}{Z2arg1})xx \= := xx\=\psdiag{3}{9}{Z2im1}
\kill
$\hat Z_f$(\pstext{Z2arg1}) \> := \> \psdiag{2}{6}{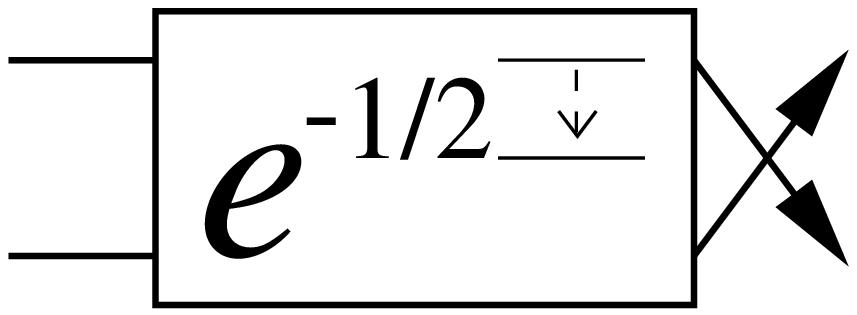}
\> $\hat Z_f$(\pstext{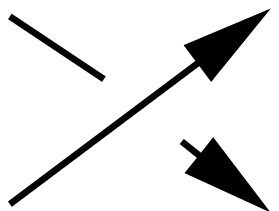}) \> := \> \psdiag{2}{6}{Z2im1} \\
[2ex]
$\hat Z_f$(\pstext{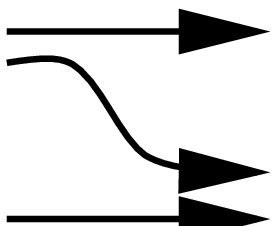}) \> := \> \psdiag{2}{6}{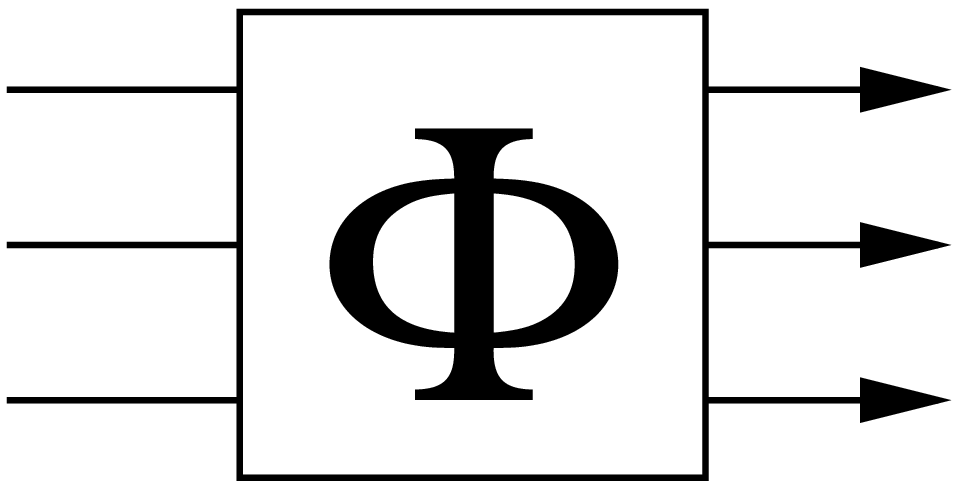}
\> $\hat Z_f$(\pstext{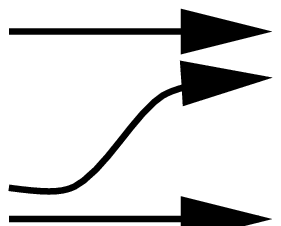}) \> := \> \psdiag{2}{6}{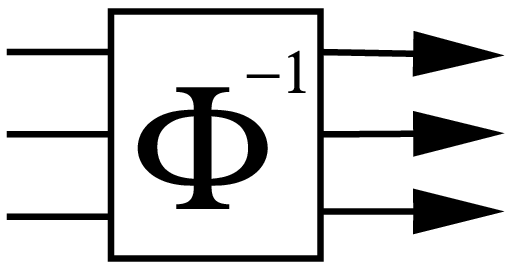} \\
[1ex]
$\hat Z_f$(\pstext{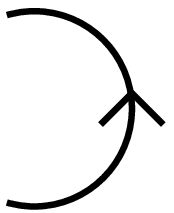}) \> := \> \psdiag{2}{6}{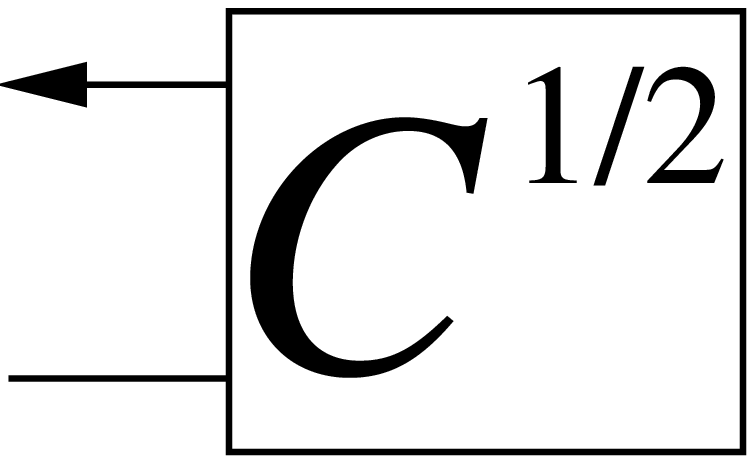}
\> $\hat Z_f$(\pstext{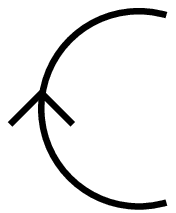}) \> := \> \psdiag{2}{6}{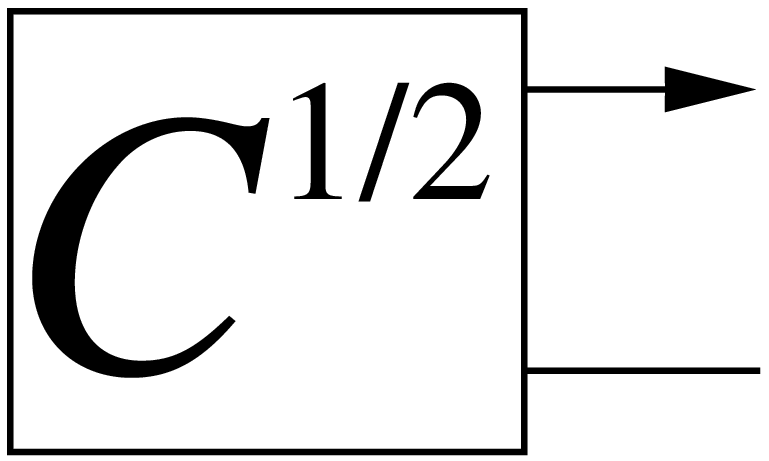} \\
[1ex]
$\hat Z_f$(\pstext{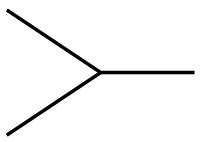}) \> := \> $\hat r_{(c_1^{x_1},c_2^{x_2},c_3^{x_3})}\cdot$\psdiag{2}{6}{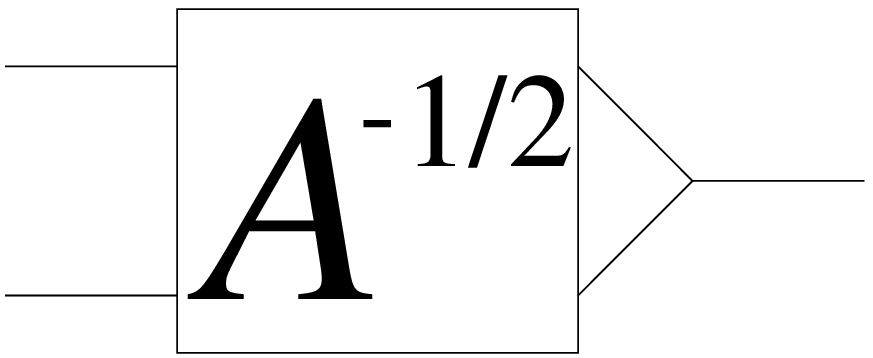}
\> $\hat Z_f$(\pstext{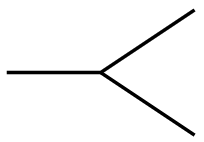}) \> := \> $\hat r_({c_1^{x_1},c_2^{x_2},c_3^{x_3}})\cdot$\psdiag{2}{6}{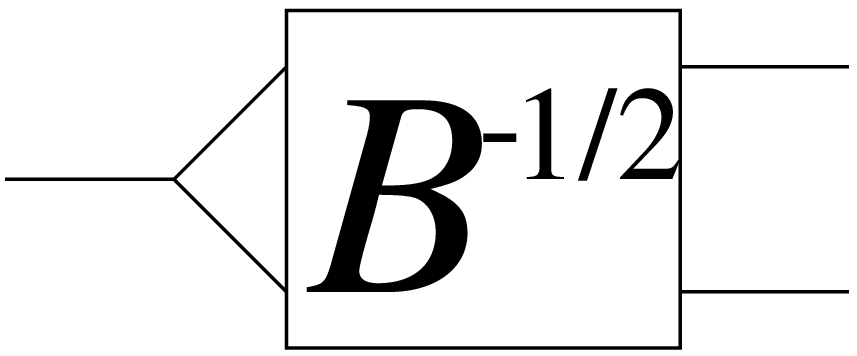} \\
\end{tabbing}
\begin{tabbing}
where xx \= \kill
where \> $e^{\pm \frac{1}{2}\psdiag{0}{3.5}{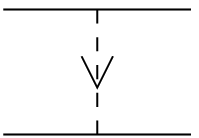}} := \sum_{n=0}^{\infty}$
$(\pm\frac{1}{2})^n\frac{1}{n!}\psdiag{1}{4}{hh2}^n$ \\[1ex]
      \> $\Phi$ is the Knizhnik-Zamolodchikov associator with arbitrary\\
      \> orientation on chords (for definition see [LM 3]) \\
      \> $C := (\pstext{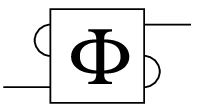})^{-1}$ \\
      \> $\hat r_{(c_1^{x_1},c_2^{x_2},c_3^{x_3})} \in \C \setminus \{ 0 \}
         \quad\forall c_1,c_2,c_3,x_1,x_2,x_3$ \\
      \> A := \psdiag{3}{10}{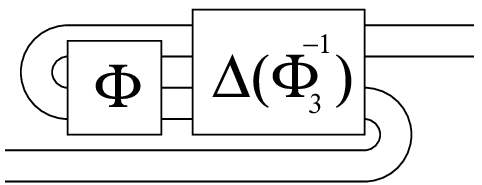}, B := \psdiag{3}{10}{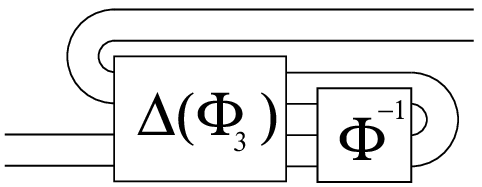}.
\end{tabbing}
For all the other generators  (i.e.\ $\pstext{Zarg1}_{u,v}$,
$\pstext{Zarg5}_{u,v}$, $\pstext{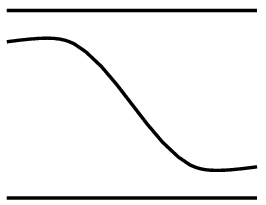}_{u,v,w}$, or
$\pstext{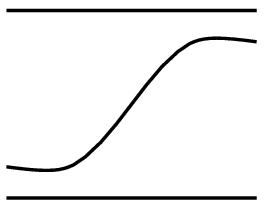}_{u,v,w}$ for some non empty non-associative coloured words u,
v, and w, and $\pstext{erzd3crb}_c,\pstext{erzd3clb}_c$ for some colour $c$) the image under $\hat Z_f$ is obtained by applying
$\Delta_u,\Delta_v$,$\Delta_w$ and $\Delta_{c^{\overline{x}}}$ to the corresponding components
of $\hat Z_f(\pstext{Z2arg1})$, $\hat
Z_f(\pstext{Z2arg5})$, $\hat
Z_f(\pstext{Z2arg2})$, $\hat
Z_f(\pstext{Z2arg6})$, $\hat Z_f(\pstext{Z2arg3})$, or
$\hat Z_f(\pstext{Z2arg7})$respectively,
e.g.\
$$\hat Z_f(\pstext{Zarg1}_{u,v}) = \Delta_v(\Delta_u(\hat Z_f(\pstext{Zarg1})_1)_{l(u)}.$$
\end{dt}

\begin{rmk}
{\em As $\Phi$, A, B, C are formal power series in certain coloured
3-diagrams with degree 0-part 1, one may take their inverses and
square roots by substituting $x$ for their higher degree parts and
expand the corresponding function of $x$ in a Taylor series. }
\end{rmk}

\begin{rmk}
{\em Since the number of trivalent vertices of a certain type in a
3-tangle is invariant, we have choices for $\hat
r_{(c_1^{x_1},c_2^{x_2},c_3^{x_3})}$. The only restriction to these
choices is, that $\hat r_{(c_1^{x_1},c_2^{x_2},c_3^{x_3})}= \hat
r_{ (c_4^{x_4},c_5^{x_5},c_6^{x_6})}$ whenever
$(c_1^{x_1},c_2^{x_2},c_3^{x_3})$ and $(c_4^{x_4},
c_5^{x_5},c_6^{x_6})$ differ by a  permutation.}
\end{rmk}

{\bf Proof of \ref{VK}:} In section 1 of [MO], Murakami and Ohtsuki
define the universal Vassiliev-Kontsevich invariant for
(uncoloured) oriented 3-tangles and prove that it is indeed an
invariant. First we modify their definition---without destroying
the invariance of the functor---in a way that keeps it well defined
for unoriented 3-nets and 3-diagrams\footnote{in the case of the
Lie algebras $\frakg_2,\frakf_4, \frake_7,$ and $\frake_8$ this
is useful because it allows us to work with unoriented 3-nets
(see remark \ref{unori})}: Omit the signs accounting for the
orientation of the strands, and introduce the antisymmetry relation
($AS$) instead of the corresponding (implicit) symmetry relation. In
order to obtain as degree 1-part of $\hat Z_f(\pstext{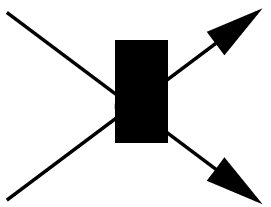}) := \hat
Z_f(\pstext{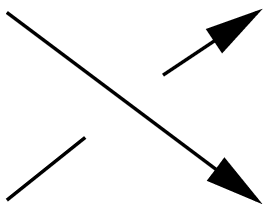} - \pstext{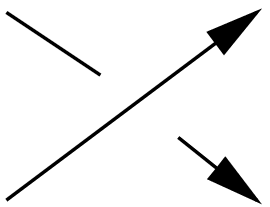})$ the diagram in which the
double point is replaced by a chord (arriving at the support on
either side like this: \pstext{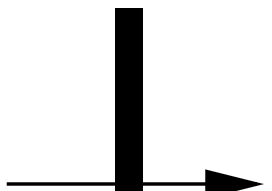}), we have adjusted the sign in
the exponent of the image of the crossings.\newline The extension
of the (uncoloured) universal Vassiliev-Kontsevich invariant to
coloured 3-nets succeeds because equivalent coloured 3-nets are equivalent
as uncoloured 3-nets and the colouring is preserved by the universal Vassiliev-
Kontsevich invariant.\schluss\\

Now we focus on the Lie algebraic part
of the construction of the invariants. As we aim at skein relations
related to the standard representation $V$ of some exceptional Lie
algebra $\frakg$, we are interested in the decomposition of
$V\otimes V$ into irreducible representations of $\frakg$. In the
following summary, definitions and notations are as
given in [H]and used in [LCL]:\\
\\
\begin{tabular}{@{} l p{12.5cm}}
$\frakg_2$:& $V:=$ the 7-dimensional irreducible representation of
$\frakg_2$\\&$V\otimes V\cong {\bf C}\oplus W\oplus V\oplus L$\\ &
with highest weights: $V$: (1,0), $L$: (0,1) (adjoint representation),
$W$: (2,0).
\end{tabular}
\\[1 ex]
\begin{tabular}{@{} l p{12.5cm}}
$\frakf_4$:& $V:=$ the 26-dimensional irreducible representation of
$\frakf_4$\\&$V\otimes V\cong {\bf C}\oplus V\oplus U\oplus L\oplus
W$\\ & with highest weights: $V$: (0,0,0,1), $L$: (1,0,0,0) (adjoint
representation), $U$: (0,0,0,2), $W$: (0,0,1,0).
\end{tabular}
\\[1ex]
\begin{tabular}{@{} l p{12.5cm}}
$\frake_6$:& $V:=$ a 27-dimensional irreducible representation of
$\frake_6$\\&$V\otimes V\cong U\oplus V^*\oplus W$\\ & with highest
weights: $V$: (0,0,0,0,0,1), $V^*$: (1,0,0,0,0,0), $U$: (0,0,0,0,0,2),
$W$: (0,0,0,0,1,0).
\end{tabular}
\\[1ex]
\begin{tabular}{@{} l p{12.5cm}}
$\frake_7$:& $V:=$ the 56-dimensional irreducible representation of
$\frake_7$\\&$V\otimes V\cong {\bf C}\oplus U\oplus L\oplus W$\\&
with highest weights: $V$: (0,0,0,0,0,0,1), $L$: (1,0,0,0,0,0,0)
(adjoint representation), $U$: (0,0,0,0,0,1,0), $W$: (0,0,0,0,0,0,2).
\end{tabular}
\\[1ex]
\begin{tabular}{@{} l p{12.5cm}}
$\frake_8$:& $V:=L$, the adjoint representation of
$\frake_8$\\&$V\otimes V\cong {\bf C}\oplus U\oplus X\oplus V\oplus
W$\\ & with highest weights: $V$: (0,0,0,0,0,0,0,1), $U$:
(1,0,0,0,0,0,0,0), $X$: (0,0,0,0,0,0,0,2), $W$: (0,0,0,0,0,0,1,0).
\end{tabular}
\\[2ex]
Note that an occurrence of ${\bf C}$ in the decomposition of
$V\otimes V$ implies that $V$ is selfdual (i.e.\ $V\cong V^*$);
so all standard representations except the one of $\frake_6$ are
selfdual.\\[2ex] We now fix for each exceptional Lie algebra
$\frakg$ a set ${\cal S}(\frakg)$ of irreducible representations of
$\frakg$:
\vspace{-5mm}
\begin{tabbing}
${\cal S}(\frakg_2):=\{{\bf C}, V,V^*,L,L^*\}$\qquad\= ${\cal
S}(\frakf_4):=\{{\bf C},V,V^*,L,L^*\}\qquad$\=${\cal
S}(\frake_6):=\{{\bf C},V,V^*,U,U^*,L,L^*\}$\kill\\ ${\cal
S}(\frakg_2):=\{{\bf C},V,V^*,L,L^*\}$\> ${\cal S}(\frakf_4):=\{{\bf
C},V,V^*,L,L^*\}$\>${\cal S}(\frake_6):=\{{\bf C},V,V^*,U,U^*,L\}$\\${\cal
S}(\frake_7):=\{{\bf C},V,V^*,L,L^*\}$\> ${\cal S}(\frake_8):=\{{\bf C},
V,V^*,U,U^*\}$
\end{tabbing}
Note that all these representations except $V,V^*,U,$ and $U^*$ in
the case $\frakg=\frake_6$ are selfdual.\newline
A check with the program
\mbox{{\sf L\hspace{-0.2em}\raisebox{0.5ex}{\scriptsize I}E}} (see [LCL])
shows that for $X,Y,Z \in {\cal S}(\frakg)$ the dimension of
Hom$_\frakg(X\otimes Y,Z)$ is 1 or 0; so we have a $\frakg$-linear
map $p^\frakg_{X\otimes Y,Z}$ from $X\otimes Y$ to $Z$ that is
either 0 or uniquely defined up to a scalar.\newline
This is important for proposition \ref{psi} and has partly
motivated the choice of the ${\cal S}(\frakg)$.

\begin{df}
The family $\{p^\frakg_{X\otimes Y,Z}:X\otimes Y\rightarrow
Z|X,Y\in {\cal S}\backslash\{\C\},Z\in{\cal S}(\frakg)\}$ of $\frakg$-linear maps is {\em
consistent} if it has the following properties:
\begin{itemize}
\item[(i)] $\pg{X}{Y}{Z}\not\equiv 0$ whenever dim Hom$_\frakg(X\otimes Y,Z)\neq
0$ and $\pg{L}{L}{{\bf C}}=h\kappa$ where $h$ is any non-zero
complex number and $\kappa$ the Killing form,
\item[(ii)] $\pg{X}{X^*}{{\bf C}}(x\otimes y)^=\left\{
\begin{array}{l}
-\pg{X^*}{X}{{\bf C}}(y\otimes x)\quad
(\forall x\in X, y\in X^*) \mbox{ if } \frakg=\frake_7 \mbox{ and } X=V\\
\pg{X^*}{X}{{\bf C}}(y\otimes x)\quad
(\forall x\in X, y\in X^*) \mbox{ otherwise,}
\end{array}\right.$
\item[(iii)] For $X,Y,Z\in {\cal S}(\frakg)\backslash\{\C\}$ $\pg{Y}{Z^*}{X^*}=(id_{X^*}\otimes \pg{Z}{Z^*}{{\bf C}})\circ
(id_{X^*}\otimes \pg{X}{Y}{Z}\otimes id_{Z^*})\circ(\ig{{\bf
C}}{X^*}{X}\otimes id_Y\otimes id_{Z^*})$ where $\ig{X}{Y}{Z}\in$ Hom$(X,Y\otimes Z)$ is the
unique $\frakg$-linear map with $\pg{Y}{Z}{X}\circ\ig{X}{Y}{Z}=id_X$ if
$\pg{Y}{Z}{X}\not\equiv 0$, and $\ig{X}{Y}{Z}\equiv 0$ otherwise.
\end{itemize}
\end{df}
Of course, there exist such consistent families.\newline

For the construction of $\Psi_\frakg$ we will also need the
following $\frakg$-linear maps:$$flip_{X\otimes Y}:X\otimes Y
\rightarrow Y\otimes X\quad (\forall X,Y \in {\cal S}(\frakg)),$$
the $\frakg$-linear map taking $x\otimes y$ to $y\otimes x\quad
(\forall x\in X,y\in Y)$.\newline

Let ${\hat{h}}$ be a formal parameter.
\begin{df}
The category ${\cal C}(\frakg)$ is the monoidal ${\bf
C}[\hspace{-0.1em}[\hat {h}]\hspace{-0.1em}]$-category with objects
{\em Obj(${\cal C}(\frakg)$)}$:=\{{\bf C}[\hspace{-0.1em}[\hat
{h}]\hspace{-0.1em}]
\otimes_\C U\,|\, U \mbox{ is a tensor product over }\C \mbox{ with  factors in }
{\cal S}(\frakg)\}$ and with the following morphism spaces:\\
\centerline{$\mbox{{\em Mor}}_{{\cal C}(\frakg)}( {\bf
C}[\hspace{-0.13em}[\hat {h}]\hspace{-0.13em}]\otimes_\C U_1, {\bf
C}[\hspace{-0.13em}[\hat {h}]\hspace{-0.13em}]\otimes_\C U_2):=
{\bf C}[\hspace{-0.13em}[\hat
{h}]\hspace{-0.13em}]\otimes_\C\mbox{{\em
Hom}}_{\frakg}(U_1,U_2).$}
\end{df}

The definition of $\Psi_\frakg$ is contained in the following
proposition whose proof will be omitted, because it just consists
in checking straightforwardly that $\Psi_\frakg$ respects all
relations required.
\begin{prop}\label{psi}
For all consistent families $\{p^\frakg_{X\otimes Y,Z}|X,Y\in {\cal S}(\frakg)\backslash \C,Z \in
{\cal S}(\frakg)\}$, there exist $k_{X,Y,Z}\in{\bf C}$ for which we
obtain a well defined $\C[\hspace{-0.1em}[\hat
{h}]\hspace{-0.1em}]$-linear monoidal functor $\Psi_\frakg: \hat{\cal
D}_{3c}
\rightarrow {\cal C}(\frakg)$ by setting
\begin{itemize}
\item[(i)]
$\Psi_\frakg(\mbox{\it empty word}):=
\C[\hspace{-0.1em}[\hat{h}]\hspace{-0.1em}]$
\begin{tabbing}
$\Psi_\frakg(so^\ra):=
\C[\hspace{-0.1em}[\hat{h}]\hspace{-0.1em}]\otimes
V;\hspace{1cm}$\=
$\Psi_\frakg(da^\la):= \C[\hspace{-0.1em}[\hat
{h}]\hspace{-0.1em}]\otimes L;$\\[1.5ex]
$\Psi_\frakg(so^\la):=
\C[\hspace{-0.1em}[\hat {h}]\hspace{-0.1em}]\otimes V^*;$\>
$\Psi_\frakg(do^\ra):=\left\{\begin{array}{l}
\C[\hspace{-0.1em}[\hat{h}]\hspace{-0.1em}]\otimes U^* \mbox{ if }
\frakg=\frake_6 \mbox{ or } \frake_8\\
\C[\hspace{-0.1em}[\hat{h}]\hspace{-0.1em}]\otimes V^* \mbox {
otherwise;}
\end{array}
\right.$
\\[1.5ex]
$\Psi_\frakg(da^\ra):=
\C[\hspace{-0.1em}[\hat {h}]\hspace{-0.1em}]\otimes L;
$\>$
\Psi_\frakg(do^\la):=
\left\{\begin{array}{l}
\C[\hspace{-0.1em}[\hat{h}]\hspace{-0.1em}]\otimes U \mbox{ if }
\frakg=\frake_6 \mbox{ or } \frake_8\\
\C[\hspace{-0.1em}[\hat{h}]\hspace{-0.1em}]\otimes V \mbox {
otherwise.}
\end{array}
\right.$
\end{tabbing}
\item[(ii)] $\Psi_\frakg(
$\objez$\hspace{-1mm}
\psdiag{2}{6}{erzd3ckr}\hspace{-1mm}$\objze$):=
\left\{\begin{array}{l}
- 1\otimes flip_{\Psi_\frakg(c_1^{x_1})\otimes\Psi_\frakg(c_2^{x_2})}\mbox{ if }
\Psi_\frakg(c_1^{x_1}),\Psi_\frakg(c_2^{x_2})\in\{V,V^*\}\\
\hspace{4.3cm} \mbox{and }
\frakg = \frake_7\\
1\otimes flip_{\Psi_\frakg(c_1^{x_1})\otimes
\Psi_\frakg(c_2^{x_2})}
\end{array}\right.$
\item[(iii)]
$\Psi_\frakg($\objx$\hspace{-1mm}\psdiag{2}{6}{erzd3crb}:=1\otimes
p^\frakg_{\Psi_\frakg(c^{\overline{x}})\otimes
\Psi_\frakg(c^x),{\bf C}}\newline
\Psi_\frakg(\psdiag{2}{6}{erzd3clb}\hspace{-1mm}$\objx$):= \mbox{{\rm dim} }\Psi_\frakg(c^x)\otimes i^\frakg_{{\bf
C},\Psi_\frakg(c^{\overline{x}})\otimes
\Psi_\frakg(c^x)},$\\
where $\ig{{\bf C}}{X}{X^*}\in$ Hom$({\bf C},X\otimes X^*)$ is the
unique $\frakg$-linear map with $\pg{X}{X^*}{{\bf C}}\circ\ig{{\bf C}}{X}{X^*}=id_{\bf C}$ if
$\pg{X}{X^*}{{\bf C}}\not\equiv 0$, and $\ig{{\bf C}}{X}{X^*}\equiv 0$ otherwise.
\item[(iv)] $\Psi_\frakg($\objez$\hspace{-1mm}\psdiag{2}{6}{erzd3crt}c_3^{x_3}):= \left\{
\begin{array}{l}
1\otimes p^\frakg_{\Psi_\frakg(c_1^{x_1})\otimes
\Psi_\frakg(c_2^{x_2}),\Psi_\frakg(c_3^{x_3})} \mbox{ if } \Psi_\frakg(c_i^{x_i})\neq L
(i=1,2,3)\\[1ex]
\hat{h}\otimes p^\frakg_{\Psi_\frakg(c_1^{x_1})\otimes
\Psi_\frakg(c_2^{x_2}),\Psi_\frakg(c_3^{x_3})} \mbox{ otherwise}
\end{array}\right.$\\[1ex]
$\Psi_\frakg(c_1^{x_1}\psdiag{2}{6}{erzd3clt}\hspace{-1mm}$\objzd$):=\left\{
\begin{array}{l}
k^\frakg_{\Psi_\frakg(c_1^{x_1}),\Psi_\frakg(c_2^{x_2}),\Psi_\frakg(c_3^{x_3})}
\otimes
i^\frakg_{\Psi_\frakg(c_1^{x_1}),\Psi_\frakg(c_2^{x_2})\otimes
\Psi_\frakg(c_3^{x_3})}\\
\mbox{ if } \Psi_\frakg(c_i^{x_i})\neq L \,\,
(i=1,2,3)\\[1ex]
k^\frakg_{\Psi_\frakg(c_1^{x_1}),\Psi_\frakg(c_2^{x_2}),\Psi_\frakg(c_3^{x_3})}
\hat{h}\otimes
i^\frakg_{\Psi_\frakg(c_1^{x_1}),\Psi_\frakg(c_2^{x_2})\otimes
\Psi_\frakg(c_3^{x_3})}\\ \mbox{ otherwise,}
\end{array}\right.$\\
where $\ig{X}{Y}{Z}\in$ Hom$(X,Y\otimes Z)$ is the
unique $\frakg$-linear map with $\pg{Y}{Z}{X}\circ\ig{X}{Y}{Z}=id_X$ if
$\pg{Y}{Z}{X}\not\equiv 0$, and $\ig{X}{Y}{Z}\equiv 0$ otherwise.
\end{itemize}
\end{prop}

\begin{rmk} {\em The factors $k^\frakg_{X,Y,Z}$ in $(iv)$ depend on the maps $p^\frakg_{X\otimes Y,Z}:X\otimes Y
\rightarrow Z$. The value of $k^\frakg_{X,Y,Z}$ for a fixed map $p^\frakg_{X\otimes Y,Z}$ can be found by
solving the equation $(p^\frakg_{Y^*\otimes Y,{\bf C}}\otimes
id_Z)=p^\frakg_{Y^*\otimes X,Z}\circ(id_{Y^*}\otimes i^\frakg_{X,Y\otimes
Z})$ representing the fact that
$\Psi_\frakg(\pstext{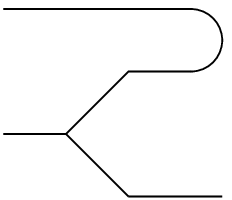})
=\Psi_\frakg(\pstext{tril})$ must hold.\newline The factor dim $\Psi_\frakg(c^x)$ in $(iii)$ has
been chosen to assure $\Psi_\frakg(\pstext{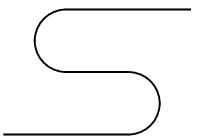})=\Psi_\frakg(\pstext{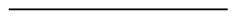})$; it
is independent of the maps
$p^\frakg_{c^{\overline{x}}\otimes c^x,{\bf C}}$.\newline The formal parameter
$\hat h$ has been introduced to assure the existence of
$\Psi_\frakg(D)$ for every morphism $D$ of $\hat {\cal D}_{3c}$: Due
to $\hat h$, the morphism spaces of ${\cal C}(\frakg)$ can be
regarded as graded spaces in the obvious way. This makes $\Psi_\frakg$ a
grade preserving functor that is well defined in every degree and
hence on the whole of $\hat {\cal D}_{3c}$.}
\end{rmk}
\begin{rmk}\label{unori}{\rm  For $\frakg=\frakg_2,\frakf_4, \frake_7$, and
$\frake_8$  the selfduality of all the representations in ${\cal{S}}(\frakg)$
allows us to identify $V,U,$ and $L$ with $V^*$,$U^*$, and $L^*$ respectively
(by fixed isomorphisms). Therefore, for any $D \in {\cal D}_{3c}\,$, $\Psi_\frakg(D)$
does not depend on the orientation of the edges and $I_{(\frakg_2,V)}$,
$I_{(\frakf_4,V)}$,$I_{(\frake_7,V)}$, and $I_{(\frake_8,V)}$ are well-defined on
unoriented 3-tangles.

}\end{rmk}

\section{Eigenvalue tables}

What we have achieved so far, is the construction of an invariant
for coloured 3-tangles---$\I:=\Psi\circ \hat Z_f$---for any exceptional
Lie algebra $\mathfrak g$ and its standard representation $V$.
Unfortunately, we cannot evaluate these invariants directly because the
expression known for the associator $\Phi$ is not explicit enough
to allow concrete computations. But we will derive skein relations
for each invariant coming from an exceptional simple Lie algebra and
its standard representation $V$. These skein relations imply recursive rules
by which we can reduce the problem of computing the invariants
to finding their values for planar coloured 3-tangles.

The idea is to cut out a small neighbourhood of a crossing and
insert something else without changing the value of the invariant.
The substitute for the crossing has to be a linear combination of
small, simple coloured 3-tangles with four univalent vertices.
Obvious candidates for such are the inverse crossing,
\mbox{\pstext{owbog}}, and, in the case of unoriented 3-nets,
\pstext{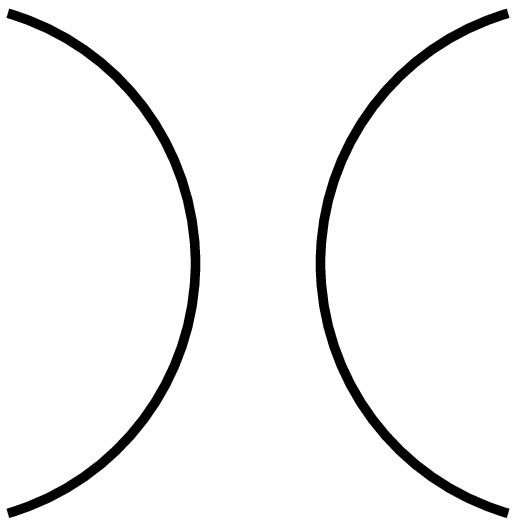}; as their values under $\I$ will prove to
be linearly independent in the endomorphism space
End$_{\mathfrak g}V\otimes V$, we include
\pstext{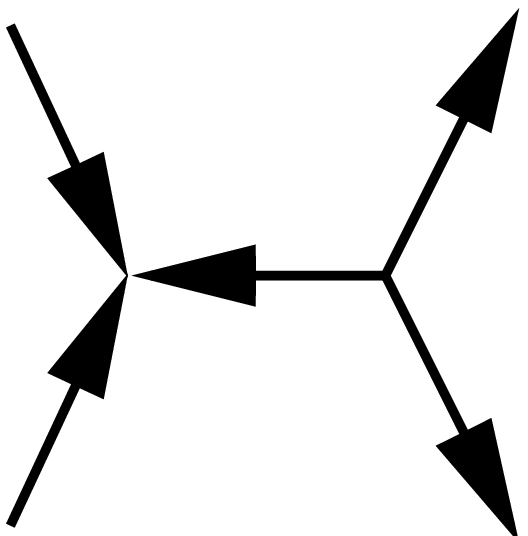} with different colourings of the edge in the middle
into our considerations.

For each $V$ considered, $V\otimes V$ decomposes into several
pairwise non-isomorphic irreducible subrepresentations. Hence
according to the lemma of Schur, the irreducible subrepresentations are
eigenspaces for each $\mathfrak g$-endomorphism
$\phi$ of $V\otimes V$, i.e.\ $\phi$ is determined by the
corresponding eigenvalues. These we must know to establish the
skein relations.

The eigenvalues of $\I (\pstext{opkreuz})$ and $\I (\pstext{onkreuz})$ are the
products of the corresponding eigenvalues of
$\Psi(e^{\mp\frac{1}{2}\psdiag{0}{2}{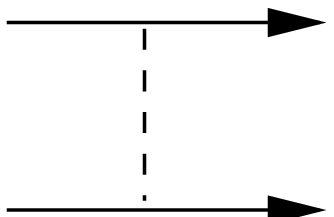}})$ and
$\Psi(\pstext{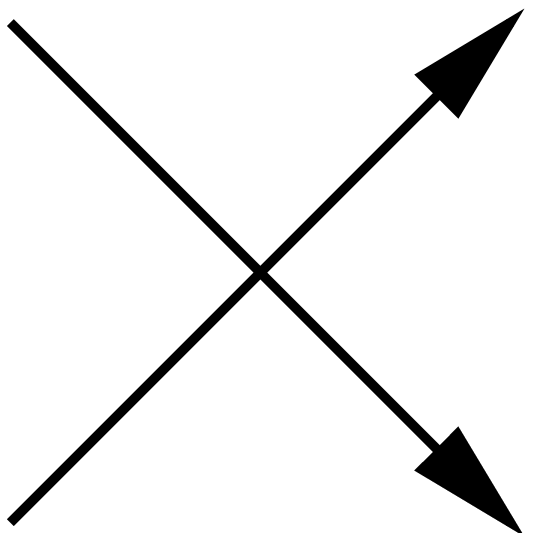})$. Therefore, we need to know the
eigenvalues of $\Psi (\pstext{ohh})$. To ascertain them we
need to know the Lie algebra and its representation explicitly.
As the corresponding tables have not been available, here a sketch of the way
how we have computed them (not mentioning the various
obstacles you meet there...):

From the Dynkin diagram characterizing the Lie algebra, the set of
roots of the Lie algebra can be derived; this is done in [SK]. Then
one can successively determine the Lie bracket of any two basis
vectors of the Lie algebra. With the help of the program
\mbox{{\sf L\hspace{-0.2em}\raisebox{0.5ex}{\scriptsize I}E}} (see [LCL]),
one can find out the heighest weight of the standard
representation $V$ of the Lie algebra in question. Letting operate
the Weyl group on this heighest weight, one gets the convex hull
of the weight lattice of the representation, hence all weights.
With the dimension of $V$ in mind (program
\mbox{{\sf L\hspace{-0.2em}\raisebox{0.5ex}{\scriptsize I}E}} !), one can
easily write down a basis of the weight spaces of $V$, whereas the
construction of a table for the operation of the Lie algebra on
$V$ is a little trickier. The last thing you need for the
description of the map $\Psi (\pstext{ohh})$ is the Killing form.

In view of the dimensions of the objects involved, it is clear
that these computations are to be made by computer. Our programs
and the various tables they produce are available at {\tt
www.math-stat.unibe.ch/$\sim$bergerab}.

The map $\I (\pstext{owbog})$ is the identity on $V\otimes V$,
hence its only eigenvalue is 1.

The eigenvalue of $\I (\pstext{sbog})$ on {\bf C} (if {\bf C} does
appear in the decomposition of $V\otimes V$) equals $\I
(\pstext{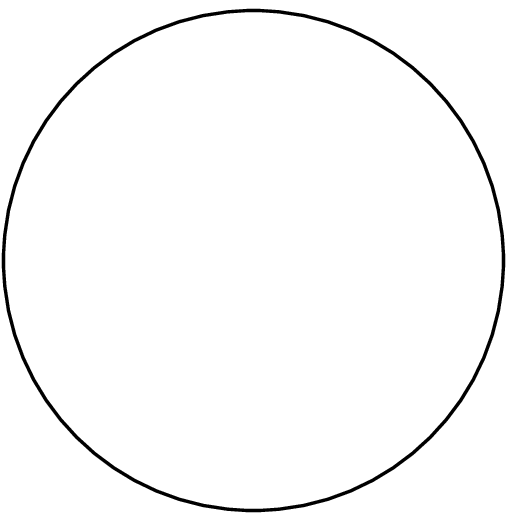})$, which is given by a formula of Rosso and Jones
(see [RJ]). On all other subrepresentations, $\I (\pstext{sbog})$
is $\equiv 0$. Let $c$ be the eigenvalue of $(\pstext{defC})^{-1}$
on ${\bf C}$. Then $\I(\pstext{kreis}:= \mbox{dim }V\cdot c$.

The value of the invariant $\I$ on a tangle of the form
\pstext{Hor} is a multiple of the projection onto the
representation corresponding to the colouring of the edge in the
middle composed with its re-imbedding. As the number of trivalent
vertices of the kind appearing in this tangle is constant under
isotopy, the factor (i.e.\ the eigenvalue on the corresponding
representation) may be chosen freely.

The following tables summarize the results of our computations for
all exceptional Lie algebras:
\vspace{4mm}\\
{\bf For the Lie algebra $\mathfrak g _2$ and its 7-dimensional
standard representation $V$} ($q=e^{-\frac{\hat h^2}{24h}}$; heighest weights of the
representations: $V$ (1,0), $L$ (0,1), $W$ (2,0)):
\vspace{6mm}\\
\centerline{\begin{tabular}{|c|c||c|c|c|c|c|c|} \hline
\multicolumn{2}{|c||}{\raisebox{0.5ex}{$\begin{array}{c}
\mbox{{\footnotesize Eigen-}} \\
\mbox{{\footnotesize value on}}\end{array}$}}
& \raisebox{0.5ex}{$\Ig$(\pstext{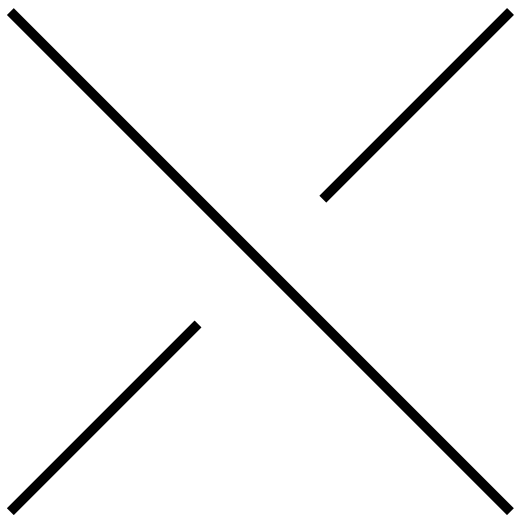})}
& \raisebox{0.5ex}{$\Ig$(\pstext{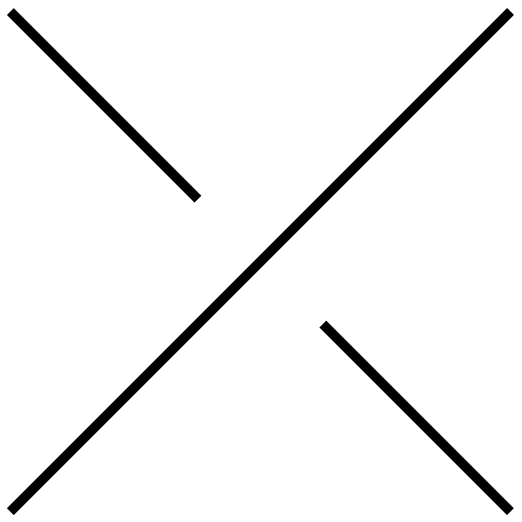})}
& \raisebox{0.5ex}{$\Ig$(\pstext{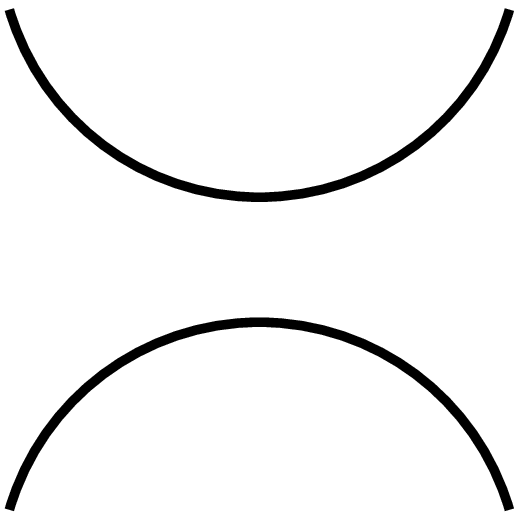})}
& \raisebox{0.5ex}{$\Ig$(\pstext{sbog})}
& \raisebox{0.5ex}{$\Ig$(\pstext{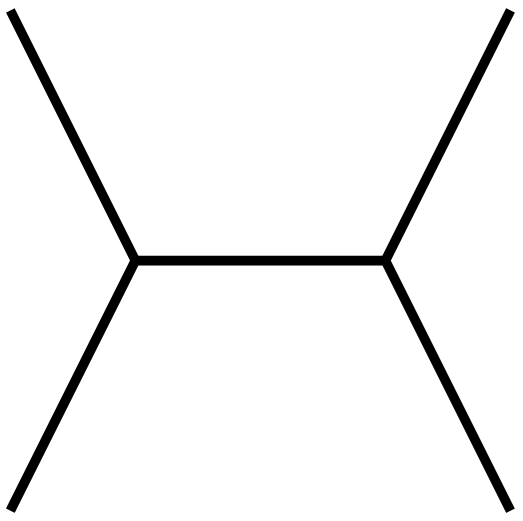})}\rule{0cm}{3.5ex}\\ \hline\hline \rule{0cm}{2.5ex}
\raisebox{-1ex}{$S^2V$} & {\bf C} & $q^{-6}$ & $q^6$ & 1 & $7c$ & 0 \\
& $W$ & $q$ & $q^{-1}$ & 1 & 0 & 0 \\ \hline \rule{0cm}{2.5ex}
\raisebox{-1ex}{$\wedge^2V$} & $V$ & $-q^{-3}$ & $-q^3$ & 1 & 0 & $r$ \\
& $L$ & $-1$ & $-1$ & 1 & 0 & 0 \\ \hline
\end{tabular}}
\vspace{6mm}\\
{\bf For the Lie algebra $\mathfrak f _4$ and its 26-dimensional
standard representation $V$} ($q=e^{-\frac{\hat h^2}{36h}}, I_{\frakf_4}:= \If$;
heighest weights of the representations:
$V$ (0,0,0,1), $L$ (1,0,0,0), $U$ (0,0,0,2), $W$ (0,0,1,0)):
\vspace{6mm}\\
\centerline{\begin{tabular}{|c|c||c|c|c|c|c|c|} \hline
\multicolumn{2}{|c||}{\raisebox{0.5ex}{$\begin{array}{c}
\mbox{{\footnotesize Eigen-}} \\
\mbox{{\footnotesize value on}}\end{array}$}}
& \raisebox{0.5ex}{$\, I_{\frakf_4}$(\pstext{pkreuz})}$\,$
& \raisebox{0.5ex}{$\, I_{\frakf_4}$(\pstext{nkreuz})}$\,$
& \raisebox{0.5ex}{$\, I_{\frakf_4}$(\pstext{wbog})}$\,$
& \raisebox{0.5ex}{$\, I_{\frakf_4}$(\pstext{sbog})}$\,$
& \raisebox{0.5ex}{$\, I_{\frakf_4}$(\pstext{H})}$\,$
& \raisebox{0.5ex}{$\quad I_{\frakf_4}$(\pstext{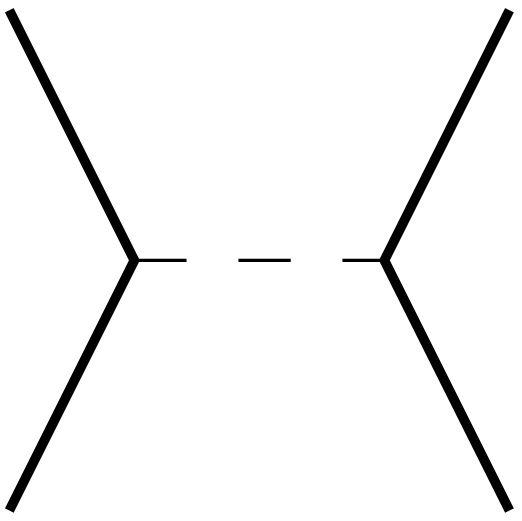})}\rule{0cm}{3.5ex}\\ \hline\hline \rule{0cm}{2.5ex}
& {\bf C} & $q^{-12}$ & $q^{12}$ & 1 & $26c$ & 0 & 0\\
$S^2V$ & $V$ & $q^{-6}$ & $q^6$ & 1 & 0 & $r$ &0\\
& $U$ & $q$ & $q^{-1}$ & 1 & 0 & 0 & 0\\ \hline \rule{0cm}{2.5ex}
\raisebox{-1ex}{$\wedge^2V$} & $L$ & $-q^{-3}$ & $-q^3$ & 1 & 0 & 0 & $s$\\
& $W$ & $-1$ & $-1$ & 1 & 0 & 0 & 0\\ \hline
\end{tabular}}
\vspace{6mm}\\
{\bf For the Lie algebra $\mathfrak e _6$ and its 27-dimensional
standard representation $V$} ($q=e^{-\frac{\hat h^2}{72h}}$; heighest weights of the
representations: $V$ (0,0,0,0,0,1), $U$ (0,0,0,0,0,2), $V^*$ (1,0,0,0,0,0),
$W$ (0,0,0,0,1,0)):
\vspace{6mm}\\
\centerline{\begin{tabular}{|c|c||c|c|c|c|c|} \hline
\multicolumn{2}{|c||}{\raisebox{0.5ex}{$\begin{array}{c}
\mbox{{\footnotesize Eigen-}} \\
\mbox{{\footnotesize value on}}\end{array}$}}
& \raisebox{0.5ex}{$\Ie$(\pstext{opkreuz})}
& \raisebox{0.5ex}{$\Ie$(\pstext{onkreuz})}
& \raisebox{0.5ex}{$\Ie$(\pstext{owbog})}
& \raisebox{0.5ex}{$\Ie$(\pstext{Hor})}
& \raisebox{0.5ex}{$\Ie$(\pstext{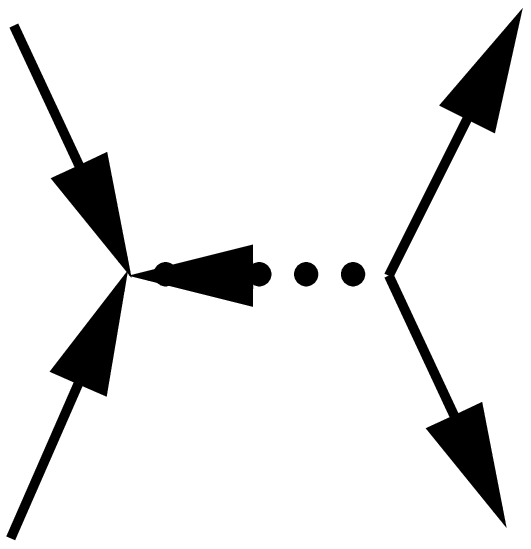})}\rule{0cm}{3.5ex}\\ \hline\hline \rule{0cm}{2.5ex}
\raisebox{-1ex}{$S^2V$} & $U$ & $q^2$ & $q^{-2}$ & 1 & 0 & $s$\\
& $V^*$ & $q^{-13}$ & $q^{13}$ & 1 & $r$ & 0 \\ \hline \rule{0cm}{2.5ex}
$\wedge^2V$ & $W$ & $-q^{-1}$ & $-q$ & 1 & 0 & 0\\ \hline
\end{tabular}
}\vspace{6mm}\\
{\bf For the Lie algebra $\mathfrak e _7$ and its 56-dimensional
standard representation $V$} ($q=e^{-\frac{\hat h^2}{144h}}$; heighest weights of the
representations: $V$ (0,0,0,0,0,0,1), $L$ (1,0,0,0,0,0,0), $U$ (0,0,0,0,0,1,0),
$W$ (0,0,0,0,0,0,2)):
\vspace{6mm}\\
\centerline{\begin{tabular}{|c|c||c|c|c|c|c|} \hline
\multicolumn{2}{|c||}{\raisebox{0.5ex}{$\begin{array}{c}
\mbox{{\footnotesize Eigen-}} \\
\mbox{{\footnotesize value on}}\end{array}$}}
& \raisebox{0.5ex}{$\Iee$(\pstext{pkreuz})}
& \raisebox{0.5ex}{$\Iee$(\pstext{nkreuz})}
& \raisebox{0.5ex}{$\Iee$(\pstext{wbog})}
& \raisebox{0.5ex}{$\Iee$(\pstext{sbog})}
& \raisebox{0.5ex}{$\Iee$(\pstext{Hstr})}\rule{0cm}{3.5ex}\\ \hline\hline \rule{0cm}{2.5ex}
\raisebox{-1ex}{$\wedge^2V$} & {\bf C} & $q^{-57}$ & $q^{57}$ & 1 & $56c$ & 0 \\
& $U$ & $q^{-1}$ & $q$ & 1 & 0 & 0 \\ \hline \rule{0cm}{2.5ex}
\raisebox{-1ex}{$S^2V$}& $L$ & $-q^{-21}$ & $-q^{21}$ & 1 & 0 & $s$\\
& $W$ & $-q^3$ & $-q^{-3}$ & 1 & 0 & 0 \\ \hline
\end{tabular}}
\vspace{6mm}\\
{\bf For the Lie algebra $\mathfrak e _8$ and its adjoint
representation $L$} as standard representation ($q=e^{-\frac{\hat h^2}{30h}}, I_{\frake_8}:=\Ieee$;
heighest weights of the representations:
$X$ (0,0,0,0,0,0,0,2),
$L$ (0,0,0,0,0,0,0,1), $U$ (1,0,0,0,0,0,0,0), $W$ (0,0,0,0,0,0,1,0)):
\vspace{6mm}\\
\centerline{\begin{tabular}{|c|c||c|c|c|c|c|c|} \hline
\multicolumn{2}{|c||}{\raisebox{0.5ex}{$\begin{array}{c}
\mbox{{\footnotesize Eigen-}} \\
\mbox{{\footnotesize value on}}\end{array}$}}
& \raisebox{0.5ex}{$I_{\frake_8}$(\pstext{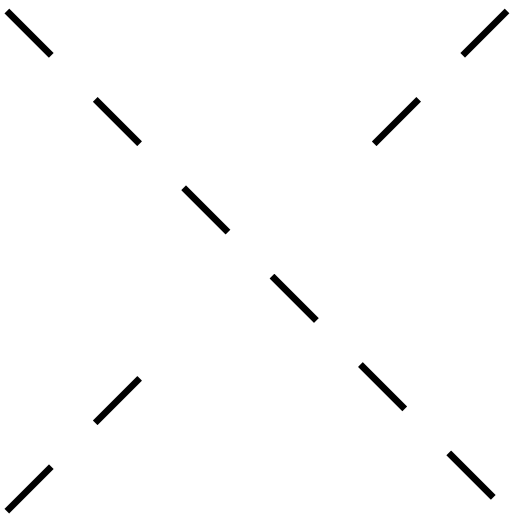})}
& \raisebox{0.5ex}{$I_{\frake_8}$(\pstext{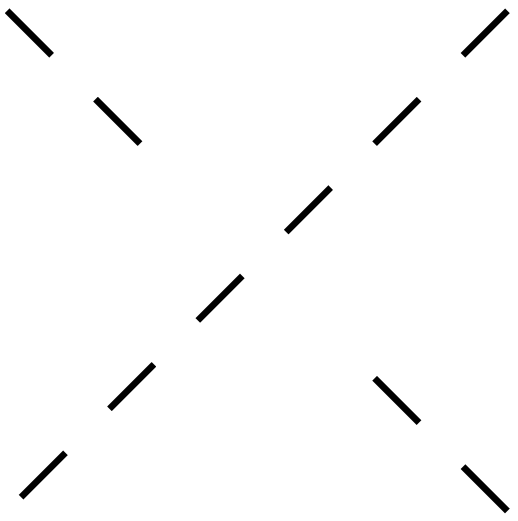})}
& \raisebox{0.5ex}{$I_{\frake_8}$(\pstext{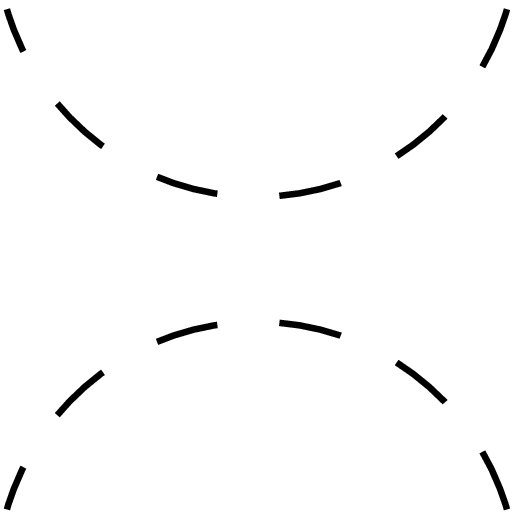})}
& \raisebox{0.5ex}{$I_{\frake_8}$(\pstext{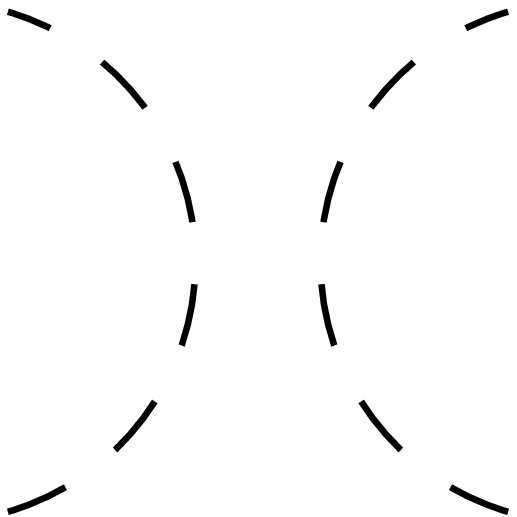})}
& \raisebox{0.5ex}{$I_{\frake_8}$(\pstext{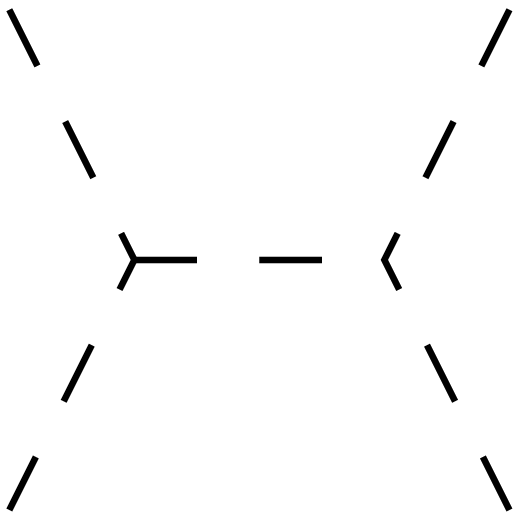})}
& \raisebox{0.5ex}{$I_{\frake_8}$(\pstext{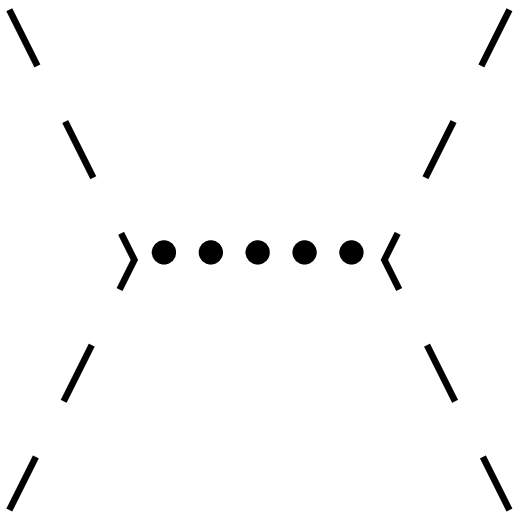})}\rule{0cm}{3.5ex}\\ \hline\hline \rule{0cm}{2.5ex}
& {\bf C} & $q^{30}$ & $q^{-30}$ & 1 & $248c$ & 0 & 0\\
$S^2L$ & $U$ & $q^{-6}$ & $q^6$ & 1 & 0 & 0 & $s$\\
& $X$ & $q$ & $q^{-1}$ & 1 & 0 & 0 & 0\\ \hline \rule{0cm}{2.5ex}
\raisebox{-1ex}{$\wedge^2L$}& $L$ & $-q^{-15}$ & $-q^{15}$ & 1 & 0 & $r$ & 0\\
& $W$ & $-1$ & $-1$ & 1 & 0 & 0 & 0\\ \hline
\end{tabular}}

\section{Skein relations}

Now it is easy to establish skein relations for the invariants
coming from the various exceptional Lie algebras: we just have to
find out how to express the eigenvalue vector corresponding to the
positive crossing as a linear combination of the eigenvalue vectors
corresponding to the other elementary tangles. But to be sure to
end up with planar coloured oriented 3-nets, skein relations
involving only one crossing are needed. In the case of the Lie
algebra $\frake_6$, we have solved this problem by introducing
\pstext{Hor}. In the remaining cases, we could of course have
proceeded by using further colours, too; but we have found rather
sophisticated methods to avoid this. We describe them in detail
for the Lie algebras $\frakg_2$ and $\frakf_4$; the cases of $\frake_7$
and $\frake_8$ can be treated in the same way as the ones of $\frakg_2$ and
$\frakf_4$, respectively. \\ \\
{\bf For the Lie algebra $\gf_2$ and its 7-dimensional irred.\
representation $V$:}
$$\Ig(\pstext{pkreuz}) = \alpha \Ig(\pstext{nkreuz})+
\beta \Ig(\pstext{wbog})
+\gamma \Ig(\pstext{sbog}) + \delta \Ig(\pstext{H})$$
where $\alpha$  :=  $q $,
      $\beta$   :=  $q-1$,
      $\gamma$  :=  $\frac{1}{7c}(-q^7+q^{-6}-q+1)$,
      $\delta$  :=  $\frac{1}{r}(q^4-q^{-3}-q+1)$.\vspace{1ex}

Since $\Ig$ is a monoidal functor and invariant under ambient isotopy,
we can deduce a new skein relation as follows:
\begin{eqnarray*}
\Ig(\pstext{nkreuz}) & = & \Ig(\pstext{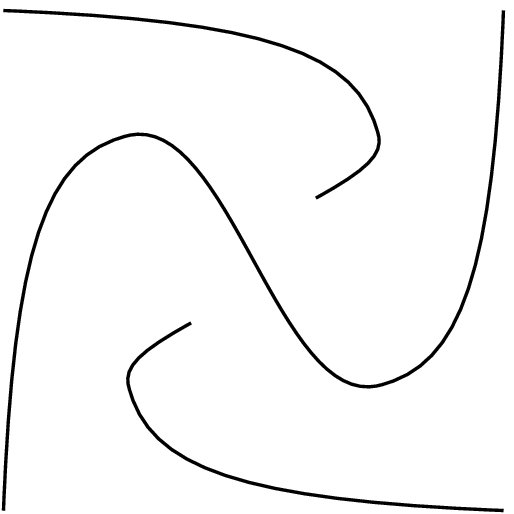}) \\
& = &  \alpha \Ig(\pstext{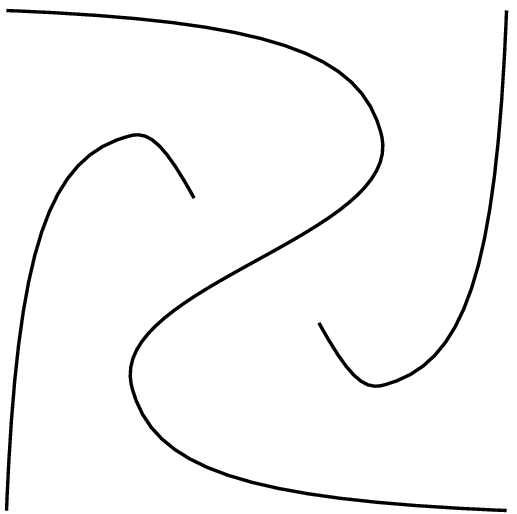})+ \beta \Ig(\pstext{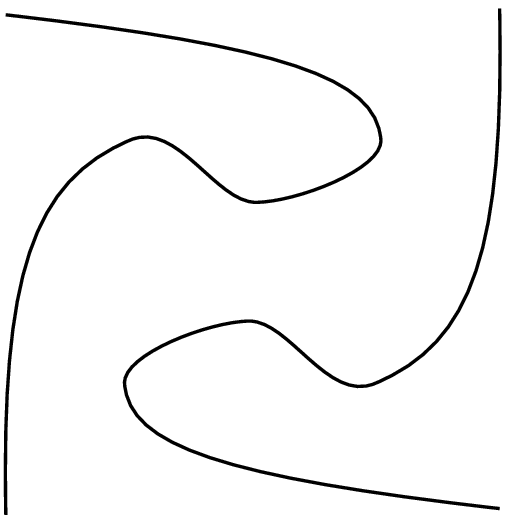})
+\gamma \Ig(\pstext{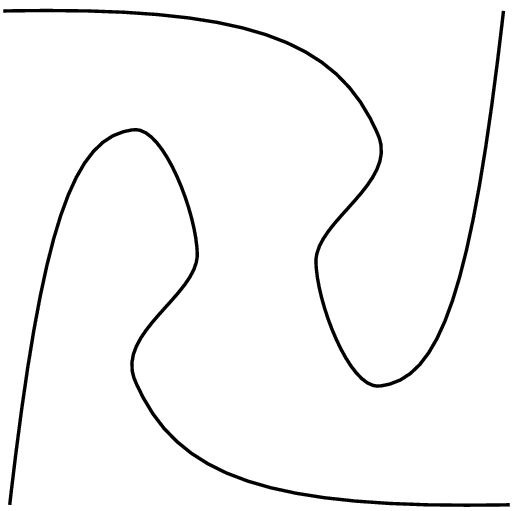}) + \delta \Ig(\pstext{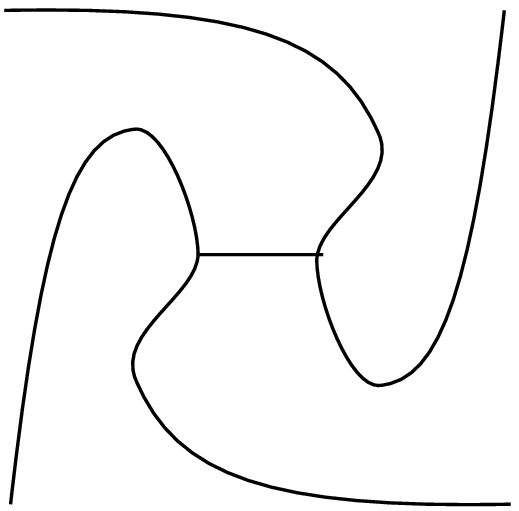}) \\
& = &  \alpha \Ig(\pstext{pkreuz})+ \beta \Ig(\pstext{sbog})
+\gamma \Ig(\pstext{wbog}) + \delta \Ig(\pstext{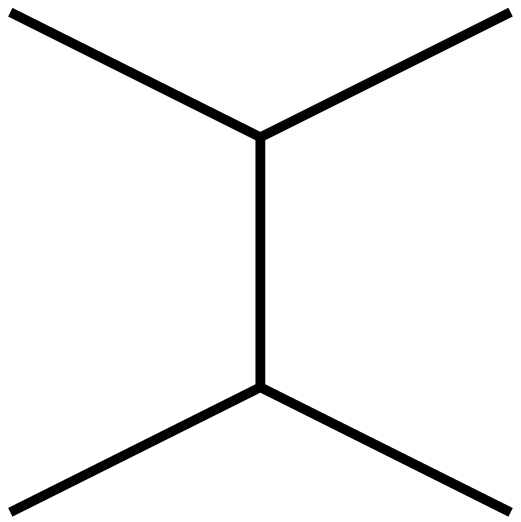}).
\end{eqnarray*}

Combining these two skein relations, we obtain:
$$\Ig(\pstext{pkreuz}) = \lambda \Ig(\pstext{wbog})+
\mu \Ig(\pstext{sbog})
+\rho \Ig(\pstext{H}) + \sigma \Ig(\pstext{I})$$
where   $\lambda$  := $\frac{q^2}{q+1}$,
        $\mu$  := $\frac{1}{q(q+1)}$,
        $\rho$  := $-\frac{q^6+q^5+q^4+q^2+q+1}{q^3(q+1)r}$,
        $\sigma$  := $q\rho.$\\ \\
{\bf For the Lie algebra $\mathfrak f_4$ and its 26-dimensional irred.\
representation $V$:}
\begin{equation}\label{SK1}
\If(\pstext{pkreuz}) = \alpha \If(\pstext{nkreuz})+
\beta \If(\pstext{wbog})
+\gamma \If(\pstext{sbog}) + \delta \If(\pstext{H})+\epsilon \If(\pstext{Hstr})
\end{equation}
where $\alpha$  :=  $q$, $\beta$ := $q-1$, $\gamma$:= $-\frac{q^9-q^8-q^5+q^4+q-1}{(q^6-q^3+1)q}$,
 $\delta$ := $-\frac{q^{13}+q^7-q^6-1}{rq^6}$, $\epsilon$ :=
 $\frac{q^7-q^4+q^3-1}{sq^3}$.\\[1ex]

To establish a skein relation involving only one crossing, the
tangles given in the table in section 2 are not enough. But we can
also  derive the eigenvalues $x_{\bf
C},x_V,x_L,x_U$, and $x_W$ for $\If(\pstext{I})$ on ${\bf C},V,L,U,$
and $W$ respectively (and as a byproduct the eigenvalues $y_{\bf C}, y_V,
y_L,y_U$, and $y_W$ for $\If(\pstext{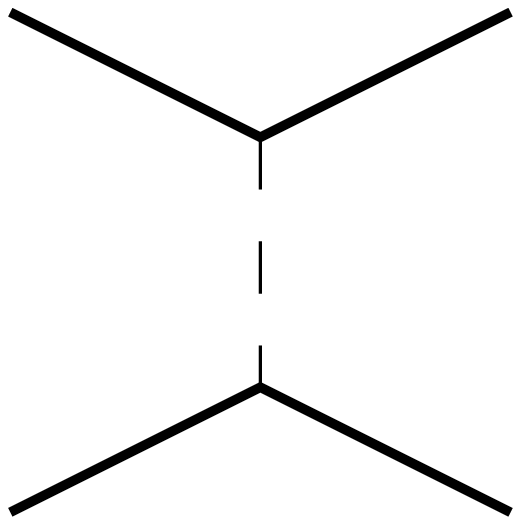}))$. To do this, we will apply the
following strategy: first we compute $x_{\bf C}$
and $y_{\bf C}$, then we use (\ref{SK1}) to establish relations
between $x_V,y_V,x_L,$ and $y_L$, and finally by examining an
ansatz for a skein relation we can set up a system of equations
whose solution leads to all the eigenvalues for $\If(\pstext{I})$
and $\If(\pstext{Istr})$.
\begin{itemize}
\item[(i)]
Eigenvalues on ${\bf C}$:\\
$x_{\bf C} = r$ (because by proceeding as in the proof of lemma
5.5 in [BS], we get

$\If(\pstext{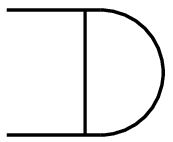})=\If(\psdiag{3}{7}{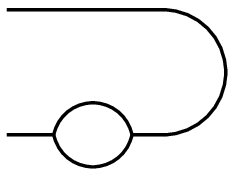})=r\If(\pstext{erzd3crb})$).\\[1.5ex]
$y_{\bf_C}=c_L^{-2}c^2s$ with $c_L =
I_{(\frakf_4,L)}(\pstext{kreis})$, because
$\If(\pstext{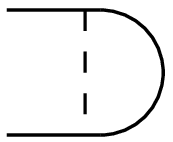})=\If(\psdiag{3}{7}{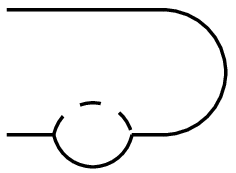})=u\If(\pstext{erzd3crb})$ and $u$
can be found by
\begin{tabbing}
$\If(\pstext{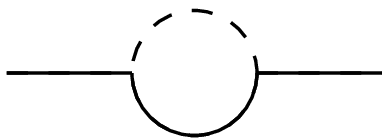})$\=$=\hat{r}^2_{(so,da,so)}\Psi_{\frakf_4}(\psdiag{3}{8}{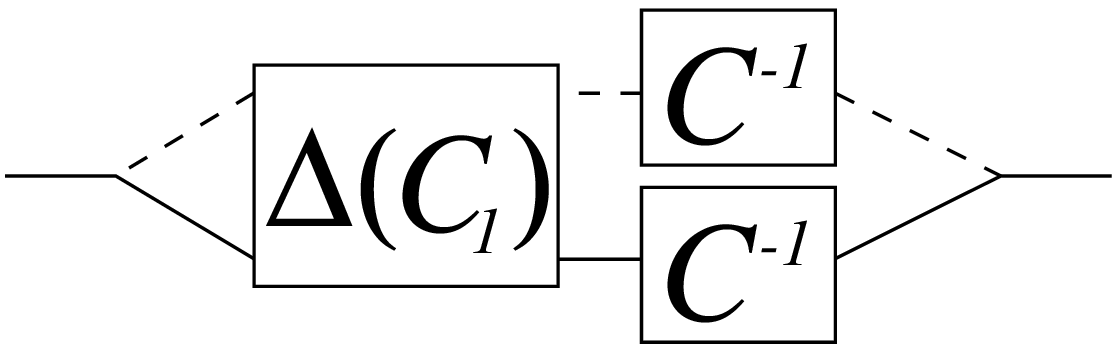})=\hat{r}^2_{(so,da,so)}
\Psi_{\frakf_4}(\psdiag{3}{8}{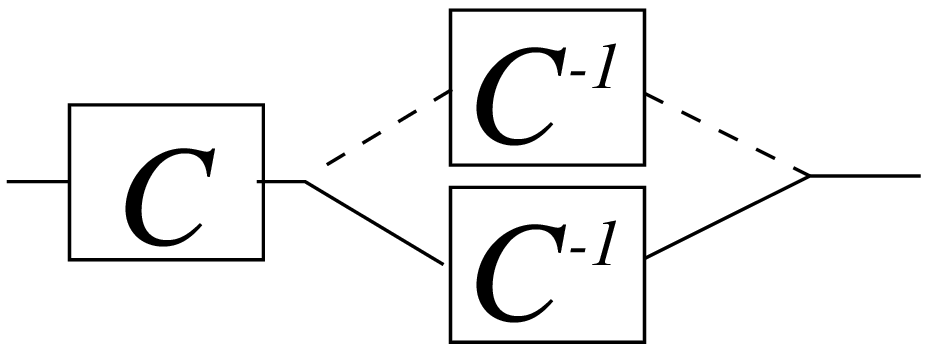})=$\\
\>$=\hat{r}^2_{(so,da,so)}c_L^{-1}\Psi_{\frakf_4}(\pstext{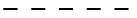})$\\
and \\
$s\Psi_{\frakf_4}(\pstext{sstrich})$\>$=\If(\pstext{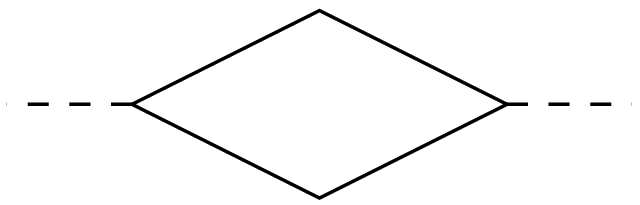})=
\hat{r}^2_{(so,da,so)}\Psi_{\frakf_4}(\psdiag{3}{8}{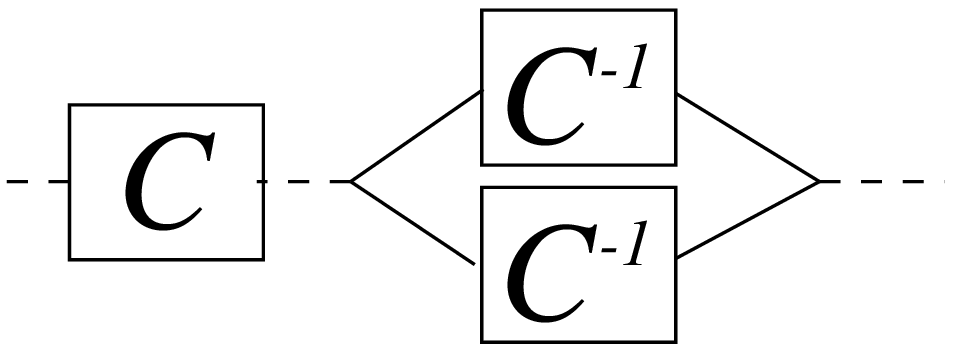})=$\\
\>$=\hat{r}^2_{(so,da,so)}c_Lc^{-2}\Psi_{\frakf_4}(\pstext{sstrich}).$
\end{tabbing}
\item[(ii)]
Eigenvalues on $V$ and $L$:\\
The following relations between the eigenvalues $x_V,y_V,x_L$ and $y_L$
hold:
$$y_V=\frac{1}{\epsilon}(q^6-\alpha q^{-6}-\gamma+\delta x_V)$$
$$y_L=\frac{1}{\epsilon}(-q^{-3}+\alpha q^3-\gamma-\delta x_L).$$
These relations have been derived by applying the skein relation
\ref{SK1}  to the crossing in \pstext{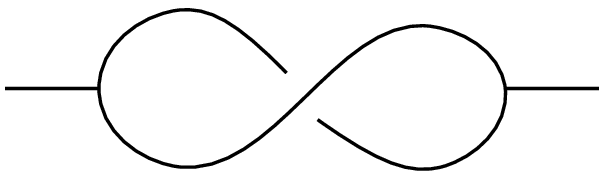} resp \pstext{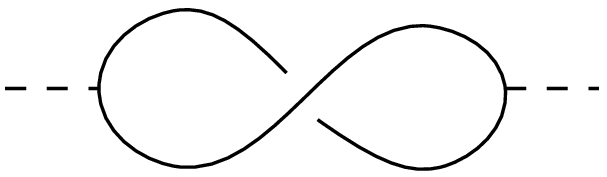}.

\item[(iii)] Note that for the skein relations
\begin{eqnarray}
\If(\pstext{Zarg1})+f_1\If(\pstext{Zarg5})
\nonumber&=& f_2\If(\pstext{wbog})+f_3\If(\pstext{sbog})+\\
&&f_4\If(\pstext{H})+f_5(\If\pstext{I})\label{SK2}
\end{eqnarray}
and\\
\vspace*{-1.34cm}
\begin{eqnarray}
\If(\pstext{Zarg1})+g_1\If(\pstext{Zarg5})
\nonumber&=& g_2\If(\pstext{wbog})+g_3\If(\pstext{sbog})+\\
& &g_4\If(\pstext{Hstr})+g_5\If(\pstext{Istr})\label{SK3}
\end{eqnarray}
turning by $90^\circ$ implies that
$f_1=1,f_2=f_3,f_4=f_5,g_1=1,g_2=g_3,$ and $g_4=g_5$.\\
By writing down (\ref{SK2}) and (\ref{SK3}) for the eigenspaces ${\bf C},V,$
 and $L$ and using the information found in $(i)$ and $(ii)$, we
get the following system of equations that can be solved by maple:
$$\begin{array}{lrcl}
{\bf C}:\hspace{1cm} &q^{-12}+q^{12}&=&f_2((1+26c)+f_4r\\
&q^{-12}+q^{12}&=&g_2(1+26c)+g_4c_L^{-2}c^2s\\
V:&q^{-6}+q^6&=&f_2(1+0)+f_4(r+x_V)\\
&q^{-6}+q^6&=&g_2(1+0)+g_4(0+\frac{1}{\epsilon}(-q^{-3}+\alpha
q^3-\gamma+\delta x_V)\\
L:&-q^{-3}-q^3&=&f_2(1+0)+f_4(0+x_L)\\
&-q^{-3}-q^3&=&g_2(1+0)+g_4(s+\frac{1}{\epsilon}(-q^{-3}+\alpha
q^3- \gamma - \delta x_L)
\end{array}$$
\end{itemize}
Knowing $f_2$ and $f_4$, we can use (2) to get all the eigenvalues
for $\If(\pstext{I})$, and once they are known, one can establish a
skein relation of the following form:
$$\If(\pstext{pkreuz}) = \lambda\If(\pstext{wbog})+
\mu\If(\pstext{sbog})+\nu\If(\pstext{H})
+\rho\If(\pstext{I})+\sigma\If(\pstext{Hstr})$$ \\
{\bf For the Lie algebra $\mathfrak e_6$ and its 27-dimensional irred.\
representation $V$:}
$$\Ie(\pstext{opkreuz}) = \alpha \Ie(\pstext{onkreuz})+
\beta \Ie(\pstext{owbog})
+\gamma \If(\pstext{Hor})$$
where $\alpha  :=  q ,
      \beta   :=  \frac{q^3-1}{q},
      \gamma  :=  -\frac{q^{27}+q^{15}-q^{12}-1}{rq^{13}}.$

$$\Ie(\pstext{opkreuz}) = \lambda \Ie(\pstext{owbog})+
\mu \Ie(\pstext{Hor})
+\rho \Ie(\pstext{Hdotor})$$
where   $\lambda$  := $-q^{-1}$,
        $\mu$  := $\frac{q^{-1}+q^{-13}}{r}$,
        $\rho$  := $-\frac{q^{-1}+q^2}{s}.$ \\ \\
{\bf For the Lie algebra $\mathfrak e_7$ and its 56-dimensional irred.\
representation $V$:}
$$\Iee(\pstext{pkreuz}) = \alpha \Iee(\pstext{nkreuz}) +
\beta \Iee(\pstext{wbog}) + \gamma \Iee(\pstext{sbog}) +
\delta \Iee(\pstext{Hstr})$$
\begin{tabbing}
where \= $\alpha := q^2 $, $\beta := \frac{q^4-1}{q}$,
      $\delta := -\frac{1}{s} \, (q^{23}+q^3-q^{-1}-q^{-21})$,\\ [1ex]
      \> $\gamma := \frac{q^{60}-q^{58}+q^{32}-q^{30}+1-q^{-2}+q^{-28}
      -q^{-30}}{(q^{16}-q^{14}+q^{12}-q^{10}+q^8-q^6+q^4-q^2+1)(q^2+1)(q^{10}+1)}$.
\end{tabbing}
$$\Iee(\pstext{pkreuz}) = \lambda \Iee(\pstext{wbog})+
\mu \Iee(\pstext{sbog})
+\rho \Iee(\pstext{Hstr}) + \sigma \Iee(\pstext{Istr})$$
\begin{tabbing}
where \= $\lambda  := -\frac{q^{88}+q^{60}+q^{57}+q^{55}+q^{47}
         +q^{45}+q^{39}+q^{37}+q^{29}+q^{28}+q^{27}+1}{
         (q^{22}-q^{16}+q^{12}+q^{10}-q^6+1)(q^2+1)^2(q^4-q^2+1)q^{28}}$,\\ [1ex]
      \> $\mu  := -\frac{q^{88}+q^{61}+q^{60}+q^{59}+q^{51}+q^{49}+q^{43}
         +q^{41}+q^{33}+q^{31}+q^{28}+1}{
         (q^{22}-q^{16}+q^{12}+q^{10}-q^6+1)(q^2+1)^2(q^4-q^2+1)q^{30}}$,\\ [1ex]
      \> $\rho  := -\frac{(q^{40}+q^{36}+q^{32}+q^{28}+q^{24}+2q^{20}+q^{16}+q^{12}
         +q^{8}+q^4+1)}{q^{21}s}$,\\ [1ex]
      \> $\sigma  := \frac{(q^4+1)(q^4-q^2+1)(q^{32}+q^{30}-q^{26}+q^{22}+q^{20}
         +q^{12}+q^{10}-q^6+q^2+1)}{q^{19}s}$.
\end{tabbing}
{\bf For the Lie algebra $\mathfrak e_8$ and its adjoint
representation $L$:}
$$\Ieee(\pstext{pkreuzstr}) = \alpha\Ieee(\pstext{nkreuzstr}) +
\beta\Ieee(\pstext{wbogstr}) + \gamma\Ieee(\pstext{sbogstr}) +
\delta\Ieee(\pstext{Hstrstr}) + \epsilon\Ieee(\pstext{Hstrdot})$$
\begin{tabbing}
where \= $\alpha:=q$, $\beta:=q-1$,
         $\gamma:=\frac{1}{248c}\,(q^{30}-q+1-q^{-29})$,\\ [1ex]
      \> $\delta:=\frac{1}{r}\,(q^{16}-q+1-q^{-15})$,
         $\epsilon:=-\frac{1}{s}\,(q^7+q-1-q^{-6})$.
\end{tabbing}

\begin{tabbing}
$\Ieee(\pstext{pkreuzstr})$ = \=$\lambda\Ieee(\pstext{wbogstr})
+ \mu \Ieee(\pstext{sbogstr})+ \nu \Ieee(\pstext{Hstrstr})$\\
[1ex]
\>$+\,\,\rho \Ieee(\pstext{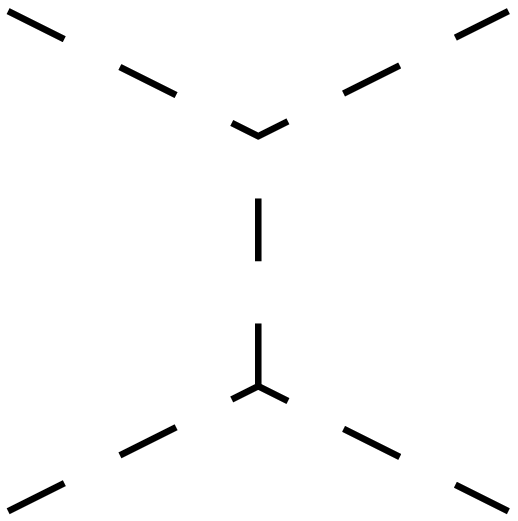})
+ \sigma \Ieee(\pstext{Hstrdot}) + \tau \Ieee(\pstext{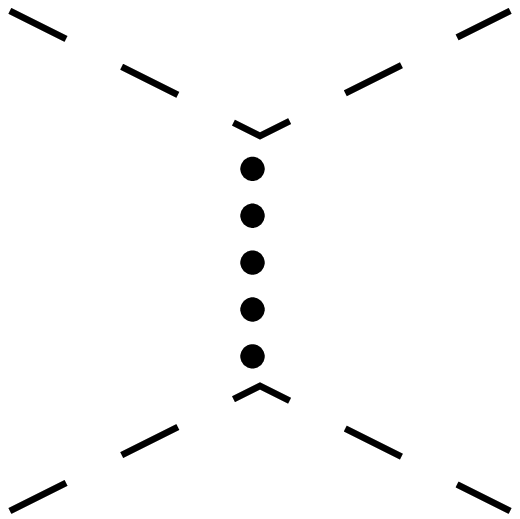}) $
\end{tabbing}

\begin{tabbing}
where \= $\lambda:=q\mu $,
         $\mu:= -\frac{1}{248c(q+1)}\sum_{k:=1}^{29}(q^k+q^{-k})-\frac{q}{c(q+1)}$, \\ [1ex]
      \> $\nu:=\frac{1}{r}\,(q^{16}-q+1-q^{-15}) $,
         $\rho:=q \nu $,
         $\sigma:=-\frac{1}{s}\,(q^7+q-1-q^{-6}) $,
         $\tau:=q\sigma $.
\end{tabbing}
Remember that $248c=\Ieee(\pstext{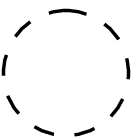})$.

\section{Reducing planar coloured 3-nets}

The idea how to reduce planar coloured 3-nets is essentially the
same as for the skein relation: cut out a small part of the 3-net
and insert something else without changing the value of the
invariant. In this case, we are to cut out a small neighbourhood
of a mesh (for definition see below), and the substitute has to be
a linear combination of planar coloured 3-nets without bounded
meshes and with the appropriate number of univalent vertices.

As $\I$ is a monoidal functor, it is enough to consider connected
coloured 3-nets; furthermore, each planar coloured 3-net is by
definition equivalent to one contained in ${\bf R}^2\times\{0\}$
with only upward pointing vectors assigned, and so we can
concentrate on this type of 3-nets.

\begin{df}
Let N be a planar coloured 3-net. A {\em mesh} of N is the closure
of a connected component of $({\bf R}^2\times\{0\})\backslash N$. A
$n-mesh$ is a mesh with n trivalent vertices in the boundary.
\end{df}

A planar coloured 3-net occuring as a result of the reduction of a
knot by means of the skein relation is closed; thus, closed,
planar coloured 3-nets are of special interest to us. About these,
we have the following lemma:

\begin{lem}([BS] section 5)
Let N be a non-empty closed, connected, planar coloured 3-net.
Then N has at least one simply connected n-mesh with $n\le 5$.
\end{lem}

Unfortunately, we are not able to resolve every possible mesh;
there may even be meshes that cannot be replaced by a linear
combination of planar coloured 3-nets without bounded meshes
(since in some cases, the dimension of the homomorphism space into
which the mesh is mapped is bigger than the number of suitable
3-nets to replace it). Nevertheless, our results make it possible
to calculate the values of these invariants on many
knots.

In the sequel, we will report the state of affairs for each case
of an exceptional simple Lie algebra.
Note that the value of $\I$ is 0 on
\objez\hspace{-1mm}\pstext{erzd3crt}$c_3^{x_3}$ if $\Psi_\frakg(c_3^{x_3})$ does
not occur as direct summand in
$\Psi_\frakg(c_1^{x_1})\otimes\Psi_\frakg(c_2^{x_2})$. Thus coloured 3-nets
containing such ``branchings'' are irrelevant for our purpose. We
will indicate for each Lie algebra the relevant ``branching
types''; with this information, it is easy to write down a
complete list of possible (relevant) $n$-meshes for any $n$.
Exemplarily, we show in the case
of $\mathfrak e_7$ how to derive the resolutions of all possible
$n$-meshes for $n\le 3$ and of a few 4-meshes.\\
\\ 
{\bf For the Lie algebra $\gf_2$ and its 7-dimensional irred.\
representation $V$:}\\
In the skein relations, only solid edges appear, and in [BS], we
have showed how to resolve every solid $n$-mesh for $n\le 5$.
Hence we can evaluate $\Ig$ on every knot.\\
\\
{\bf For the Lie algebra $\mathfrak f_4$ and its 26-dimensional irred.\
representation $V$:}\\
Relevant branching types: \pstext{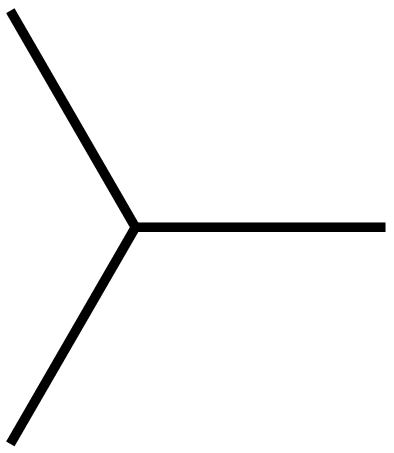},
\pstext{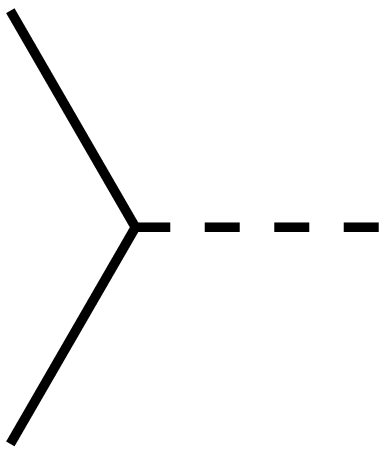}, and \pstext{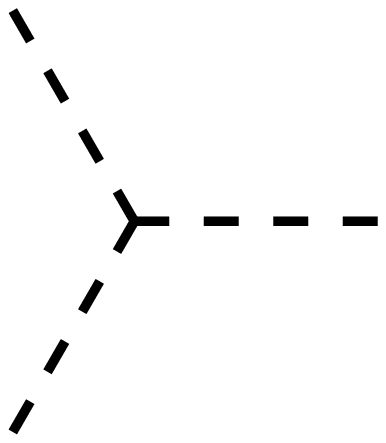}.\\
We know how to resolve any $n$-mesh for $n\le 3$ and some
4-meshes.\\
\\
{\bf For the Lie algebra $\mathfrak e_6$ and its 27-dimensional irred.\
representation $V$:}\\
Relevant branching types:
\pstext{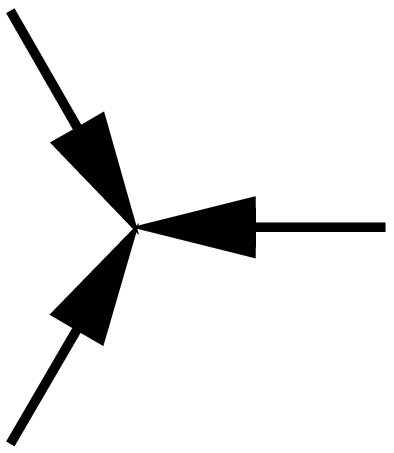}, \pstext{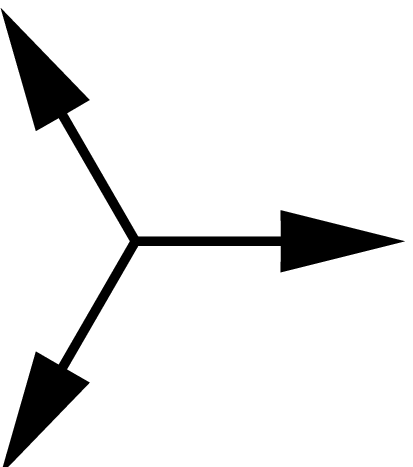},
\pstext{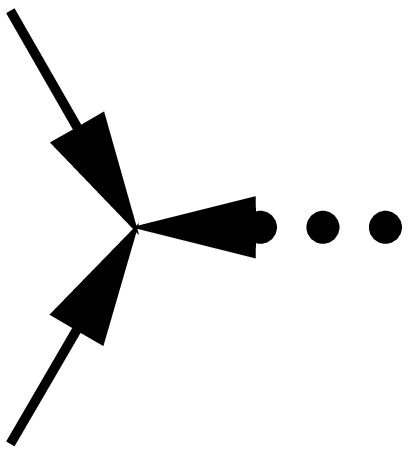}, \pstext{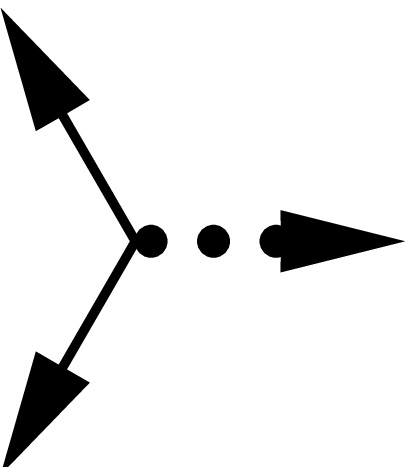}, \pstext{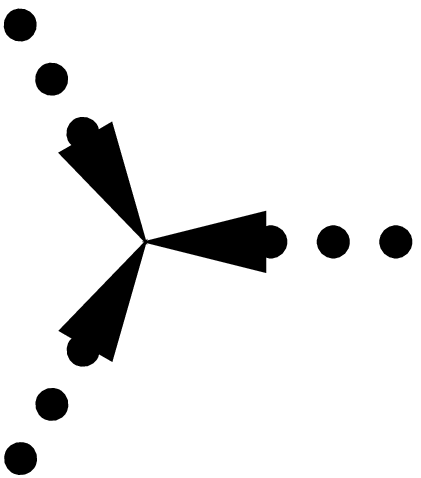}, and
\pstext{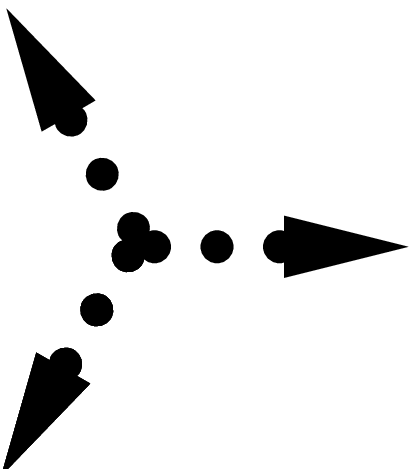}.\\
Because of the orientations, only $n$-meshes with $n$ pair can
occur.\\
We know how to resolve any 0- and 2-mesh and some 4-meshes.\\
\\
{\bf For the Lie algebra $\mathfrak e_7$ and its 56-dimensional irred.\
representation $V$:}\\
Relevant branching types: \pstext{sosoda} and \pstext{dadada}.\\
Let us denote by $u$ the eigenvalue of $\Iee(\pstext{Hstrstr})$ on
$L$, by $f$ the eigenvalue of $\Iee(\pstext{pkreuz})$ on {\bf C},
and by $l$ the eigenvalue of $\Iee(\pstext{pkreuz})$ on $L$.\\
{\bf 0-meshes:} \pstext{kreis} and \pstext{strkreis}.\\
$\Iee(\pstext{kreis})$ and $\Iee(\pstext{strkreis})$
can be computed by the formula of Rosso and Jones (see [RJ]). \\
{\bf 1-meshes:} \pstext{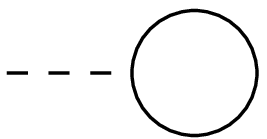} and \pstext{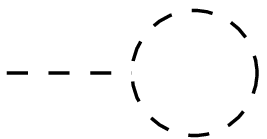}.\\
$\Iee(\pstext{sun1so})$ and $\Iee(\pstext{sun1da})$
are elements of
Hom$_{\mathfrak e_7}(L,{\bf C})$ and therefore $\equiv 0$. \\
{\bf 2-meshes:} \pstext{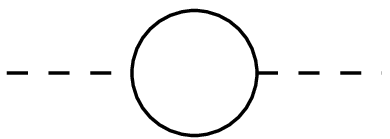}, \pstext{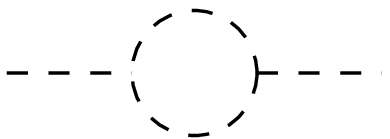}, and \pstext{fc5}.\\
$\Iee(\pstext{sun2daso})$ =
$\Iee(\pstext{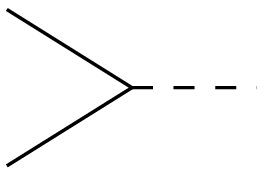})\circ\Iee(\pstext{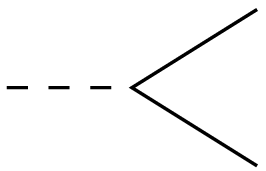})$ \\
= (eigenvalue of $\Iee(\pstext{trisosoda})$ on
$L$)(eigenvalue of $\Iee(\pstext{tridasoso})$ on
$L$)$\Iee(\pstext{sstrich})$\\
= (eigenvalue of $\Iee(\pstext{Hstr})$ on $L$)$\Iee(\pstext{sstrich})
= s \Iee(\pstext{sstrich})$.\\
Analogously: $\Iee(\pstext{sun2da}) = u \Iee(\pstext{sstrich})$.\\
By means of the skein relation, we obtain: \\
$\Iee(\pstext{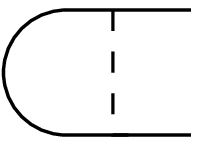})
= \frac{1}{\sigma}(\Iee(\pstext{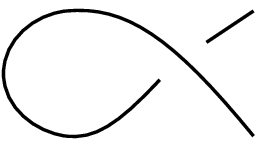})-\lambda\Iee(\pstext{erzd3clb})-
\mu c\Iee(\pstext{erzd3clb})-\rho\Iee(\pstext{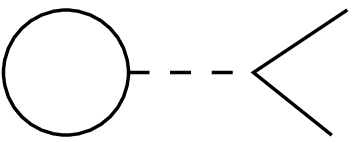})) \\
= \frac{1}{\sigma}(f-\lambda+\mu c)\Iee(\pstext{erzd3clb})$.\\
Therefore, we have:
$\Iee(\pstext{fc5}) =
\Iee(\pstext{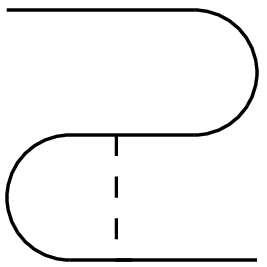}) = \frac{1}{\sigma}(f-\lambda+\mu c)\Iee(\pstext{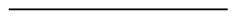})$.\\
{\bf 3-meshes:} \pstext{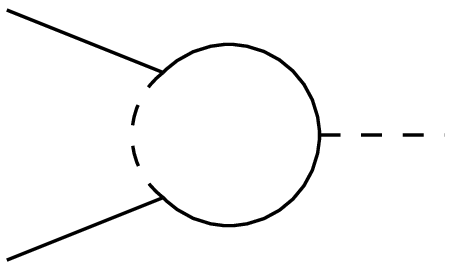}, \pstext{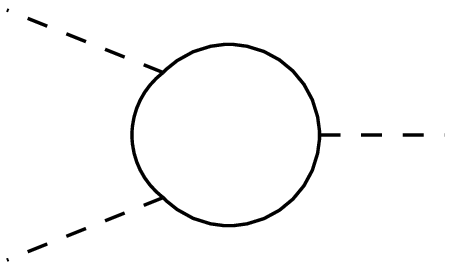}, \pstext{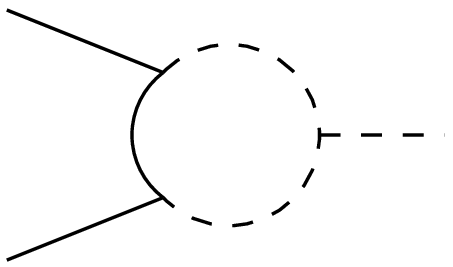}, and \pstext{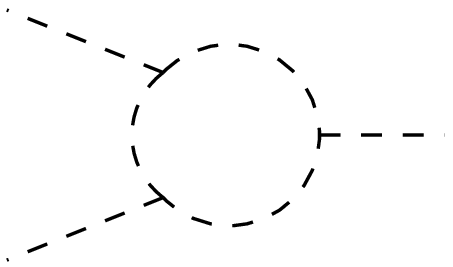}\\
The skein relation yields:
$\Iee(\pstext{sun3ssd}) = \\
\frac{1}{\sigma}(\Iee(\psdiag{1}{4}{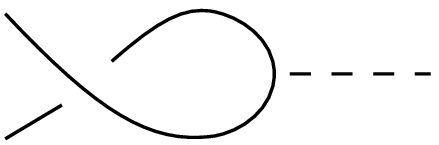})-
\lambda\Iee(\pstext{trisosoda}\!)-\mu\Iee(\pstext{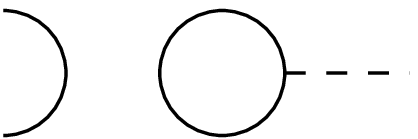})-\rho\Iee(\pstext{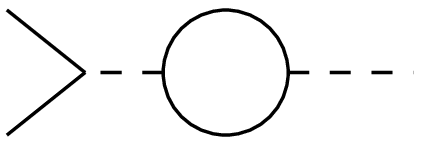}))\\
=\frac{1}{\sigma}(l-\lambda-\rho s)\Iee(\pstext{trisosoda})$.\\
The cases of \pstext{sun3dsd} and \pstext{sun3sdd} are treated
simultaneously. By setting
$\Iee(\pstext{sun3dsd})=:$\newline$t\Iee(\pstext{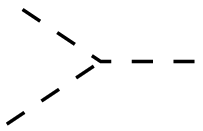})$ and
$\Iee(\pstext{sun3sdd})=:k\Iee(\pstext{trisosoda})$, we obtain
the two equations (for the second, we use the resolution of the
4-mesh \pstext{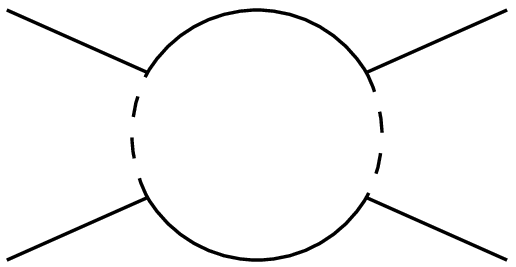} given below):\\
$ks\Iee(\pstext{sstrich})=k\Iee(\pstext{sun2daso})=\Iee(\pstext{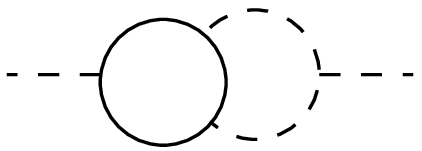})
=t\Iee(\pstext{sun2da})\\
=tu\Iee(\pstext{sstrich}) \Rightarrow ks=tu$\\
$gkt\Iee(\pstext{kreis})=kt\Iee(\pstext{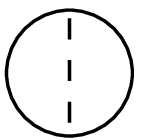})=t\Iee(\pstext{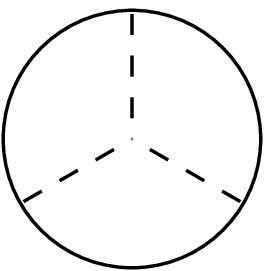})
=\Iee(\pstext{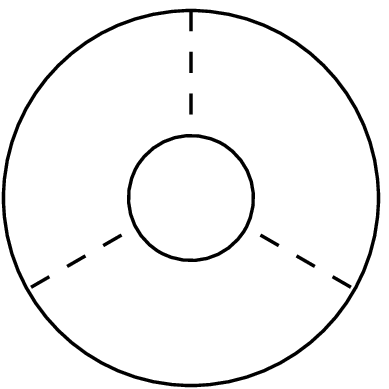})\\
=x_1\Iee(\pstext{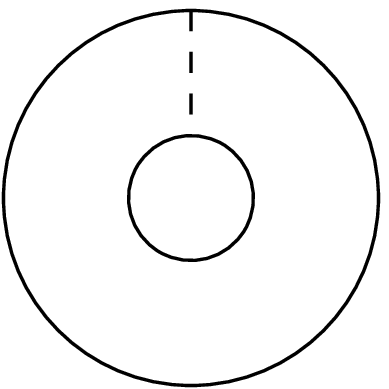})+x_2\Iee(\pstext{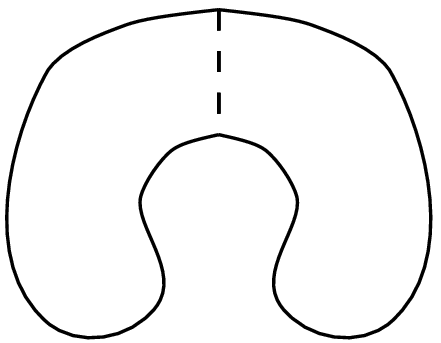})
+x_3\Iee(\pstext{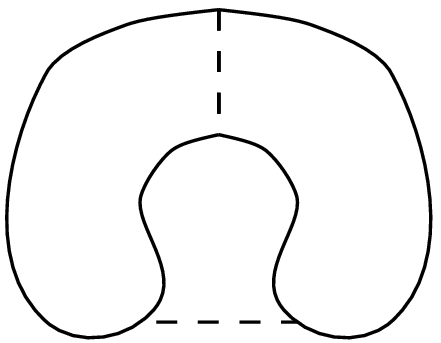})+x_4\Iee(\pstext{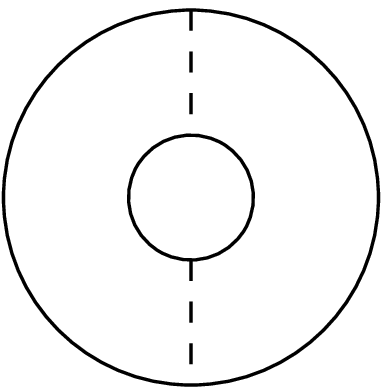})\\
=(gx_2+ghx_3+gsx_4)\Iee(\pstext{kreis}) \Rightarrow kt=x_2+hx_3+sx_4$.\\
This system of equations for $t$ and $k$ can be solved by maple.\\
The last 3-mesh, \pstext{sun3ddd}, can be resolved by means of a
skein relation involving only dashed 3-tangles; one obtains a
multiple of \pstext{tridadada}.\\
{\bf 4-meshes:} \pstext{sun4sds}, \pstext{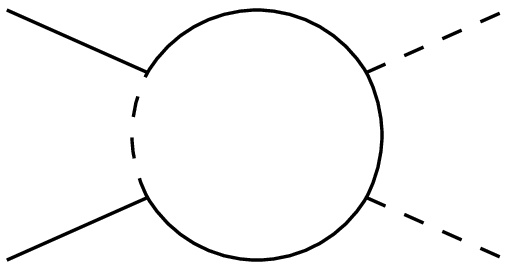}, \pstext{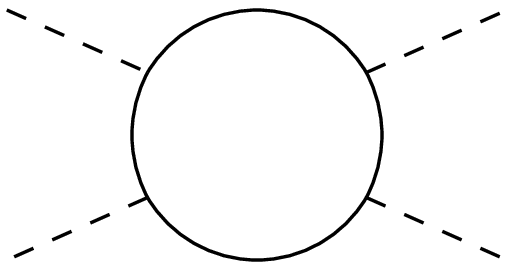},
\pstext{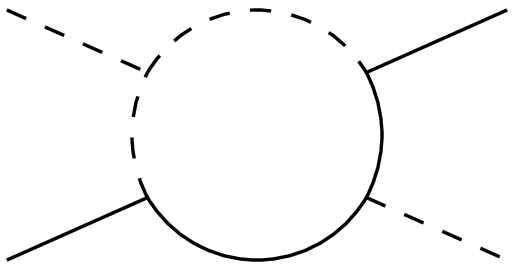}, \pstext{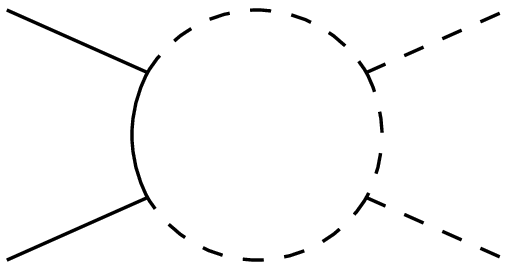}, and \pstext{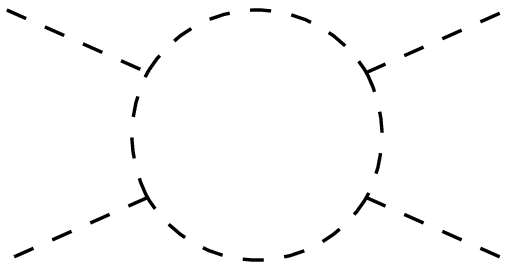}.\\
 $\Iee(\pstext{sun4sds})=x_1\Iee(\pstext{wbog})
+x_2\Iee(\pstext{sbog})+x_3\Iee(\pstext{Hstr})+x_4\Iee(\pstext{Istr})$,
where the coefficients $x_1$, $x_2$, $x_3$, and $x_4$ can be
determined by applying the skein relations.\\
The skein relation yields that\\
$ \Iee(\pstext{4ssdd0}) = \frac{1}{\sigma}(\Iee(\pstext{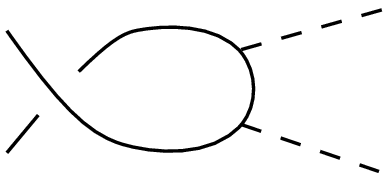})
-\lambda\Iee(\pstext{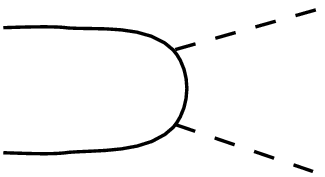})-\mu s\Iee(\pstext{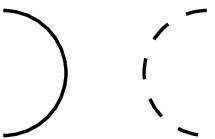})
-\rho t\Iee(\pstext{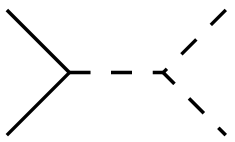}))$. \\
We will now derive a
planar substitute for \pstext{4ssdd2}. As $V\otimes V$ and
$L\otimes L$ have the common direct summands {\bf C}, $L$, and
$U$, we can make the following ansatz: \\
$\Iee(\pstext{4ssdd2}) =
X\Iee(\pstext{4ssdd4})+Y\Iee(\pstext{4ssdd5})+Z\Iee(\pstext{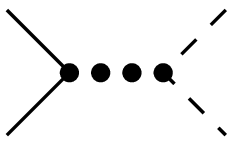})$.\\
First, we determine $X$ and $Y$:\\
$fs\Iee(\pstext{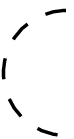})=f\Iee(\pstext{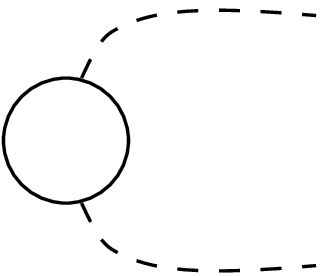})=\Iee(\pstext{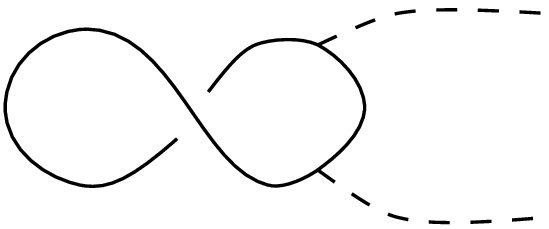})
=X\Iee(\pstext{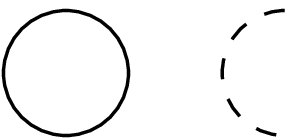}) \Rightarrow X=\frac{fs}{c}$\\
$lt\Iee(\pstext{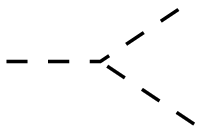})=l\Iee(\pstext{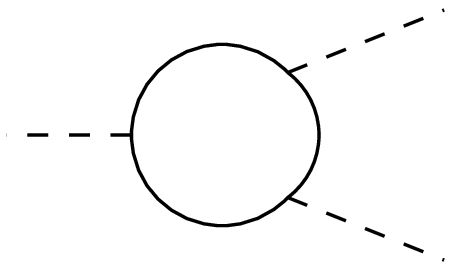})=\Iee(\pstext{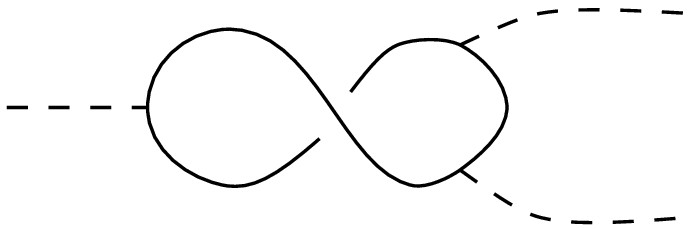})
=Y\Iee(\pstext{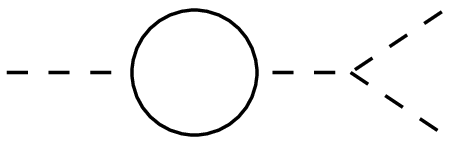})\\
=sY\Iee(\pstext{trdadada}) \Rightarrow Y=\frac{lt}{s}$.\\
By means of the following equation, we can substitute
\pstext{4ssdd6}:\\
$\Iee(\pstext{4ssdd3})=X\Iee(\pstext{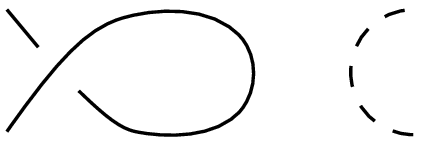})+Y\Iee(\pstext{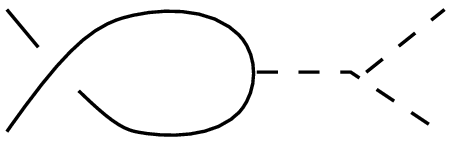})
+Z\Iee(\pstext{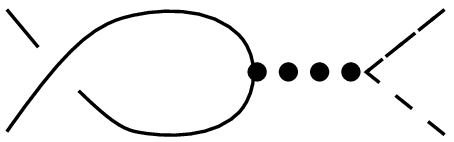})\\
=f^{-1}X\Iee(\pstext{4ssdd4})+
l^{-1}Y\Iee(\pstext{4ssdd5})+qZ\Iee(\pstext{4ssdd6})$, \\
and thereby
derive the coefficients $y_1$, $y_2$, and $y_3$ in:\\
$\Iee(\pstext{4ssdd0})=y_1\Iee(\pstext{4ssdd4})
+y_2\Iee(\pstext{4ssdd5})+y_3\Iee(\pstext{4ssdd3})$.\\
\\
{\bf For the Lie algebra $\mathfrak e_8$ and its adjoint representation $L$:}\\
Relevant branching types: \psdiag{2}{6}{dadada}, \psdiag{2}{6}{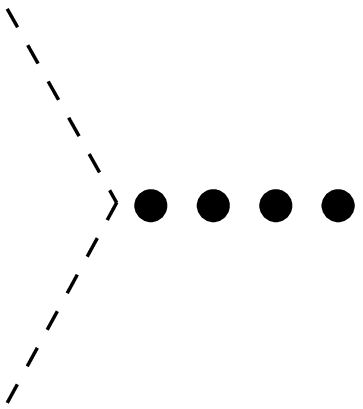},
\psdiag{2}{6}{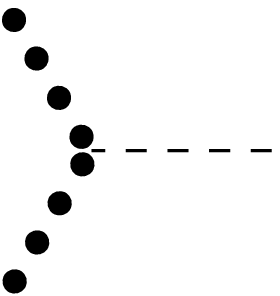} and \psdiag{2}{6}{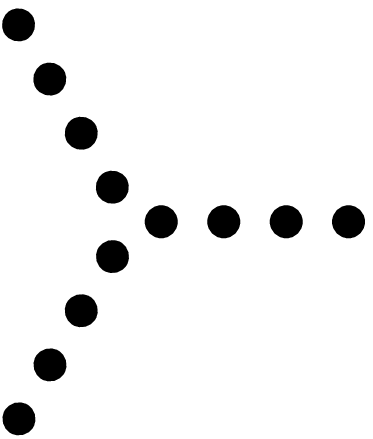}.\\
We know how to resolve all 2-meshes and those 3-meshes that do not contain
branchings of the form \psdiag{2}{6}{dododa}.

\section{References}

\begin{tabular}{@{} l p{13cm}}
[BN 1] & D. Bar-Natan, {\em On the Vassiliev knot invariants}, Topology, 34
(1995), 423-472. \\
{[}BS] & A.-B. Berger and I. Stassen, {\em The skein relation for the
$(\frakg_2,V)$-link invariant}, preprint.\\
{[}FH]  & W. Fulton and J. Harris, {\sl Representation theory}, Graduate Texts in
Mathematics \# 129, Springer-Verlag 1991.\\
{[}H]   & J.E. Humphreys, {\sl Introduction to Lie algebras and representation
theory}, Graduate Texts in Mathematics \# 9, Springer Verlag 1994.\\
{[}LCL] & M.A.A. van Leeuwen, A.M. Cohen, B. Lisser, {\em program}
\mbox{{\sf L\hspace{-0.2em}\raisebox{0.5ex}{\scriptsize I}E}},
software package for Lie group theoretical computations,
{\tt http://wallis.univ-poitiers.fr/$\sim$maavl/LiE/}.\\
{[}LM 1] & T.T.Q. Le and J. Murakami, {\em Kontsevich's integral for the
Homfly polynomial and relations between values of multiple zeta functions},
Topology Appl., 62 (1995), 193-206.\\
{[}LM 2] & T.T.Q. Le and J. Murakami, {\em Kontsevich's integral for the
Kauffman polynomial}, Nagoya Math.\ J., 142 (1996), 39-65.\\
{[}MO]  & J. Murakami and T. Ohtsuki, {\em Topological quantum field theory for
the universal quantum invariant}, Commun.\ Math.\ Phys., to appear.\\
{[}RJ]  & M. Rosso and V. Jones, {\em On the invariants of torus knots derived
from quantum groups}, J.\ knot theory ramifications, 2 (1993), 97-112.\\
{[}SK] & M. Sato and T. Kimura, {\em A classification of
irreducible prehomogeneous vector spaces and their relative
invariants}, Nagoya Math.\ J., 65 (1977), 1-155.\\
{[}V]   & P. Vogel, {\em Algebraic structures on modules of diagrams},
Invent.\ Math., to appear.
\end{tabular}

\end{document}